\documentclass[12pt]{article}
\usepackage{amsmath,amssymb,hhline}
\usepackage[russian]{babel}

\usepackage[T2A]{fontenc}    
\usepackage[cp1251]{inputenc}
\usepackage[russian]{babel} 

 \newtheorem{Th}{\quad \ Теорема}
 \newtheorem{Lemma}{\quad \ Лемма}

\begin{document}

\begin{center}
{\large\bf Р.~М.~Тригуб \qquad R. M. Trigub}
\end{center}

\begin{center}
{\large \bf О преобразовании Фурье функций двух переменных, \\ зависящих лишь от
максимума модуля этих переменных\\ \quad\\ On the Fourier transform of function of two variables
which depend only on the maximum of these variables}
\end{center}

\begin{abstract}
Для функций $f(x_{1},x_{2})=f_{0}\big(\max\{|x_{1}|,|x_{2}|\}\big)$ из
$L_{1}(\mathbb{R}^{2})$ указаны достаточные и необходимые условия принадлежности
$L_{1}(\mathbb{R}^{2})$ преобразования Фурье $\widehat{f}$ и $L_{1}(\mathbb{R}^{1}_{+})$
функции $t\cdot \sup\limits_{y_{1}^{2}+y_{2}^{2}\geq
t^{2}}\big|\widehat{f}(y_{1},y_{2})\big|$. А положительность $\widehat{f}$ на $\mathbb{R}^{2}$
полностью сведена к такому же вопросу в $\mathbb{R}^{1}$ функции
$f_{1}(x)=|x|f_{0}\big(|x|\big)+\int\limits_{|x|}^{\infty}f_{0}(t)dt$. \end{abstract}

{\bf Abstract.} For functions $f(x_{1},x_{2})=f_{0}\big(\max\{|x_{1}|,|x_{2}|\}\big)$ from
$L_{1}(\mathbb{R}^{2})$, sufficient and necessary conditions for the belonging of their Fourier transform
$\widehat{f}$ to $L_{1}(\mathbb{R}^{2})$ as well as of a function $t\cdot \sup\limits_{y_{1}^{2}+y_{2}^{2}\geq
t^{2}}\big|\widehat{f}(y_{1},y_{2})\big|$ to $L_{1}(\mathbb{R}^{1}_{+})$. As for the positivity of $\widehat{f}$ on
$\mathbb{R}^{2}$, it is completely reduced to the same question on $\mathbb{R}^{1}$ for a function
$f_{1}(x)=|x|f_{0}\big(|x|\big)+\int\limits_{|x|}^{\infty}f_{0}(t)dt$.

\textbf{Ключевые слова:} преобразование Фурье, винеровская банахова алгебра, модуль
непрерывности, целая функция, положительная определенность, экстремальные сплайны по радиальному базису


\begin{center}
    \textbf{\S0 Введение}
\end{center}

Пусть $f\in L_{1}(\mathbb{R}^{2})$ и
$f(x_{1},x_{2})=f_{0}\big(\max\{|x_{1}|,|x_{2}|\}\big)$.

В статье изучаются следующие вопросы:

1) Когда преобразование Фурье $\widehat{f}\in L_{1}(\mathbb{R}^{2})$?

2) Когда $t\cdot \sup\limits_{y_{1}^{2}+y_{2}^{2}\geq
t^{2}}\big|\widehat{f}(y_{1},y_{2})\big|\in L_{1}(\mathbb{R}^{1}_{+})$?
$\qquad\mathbb{R}^{1}_{+}=\mathbb{R}_{+}=[0,+\infty)$

3) Когда $\widehat{f}(y_{1},y_{2})\geq0$ или $>0$ всюду на $\mathbb{R}^{2}$?

По первому вопросу найдены необходимые и отдельно достаточные условия, которые для
выпуклых на отрезке функций $f_{0}$ совпадают.

Ответ по второму вопросу оказался отрицательным ($f\equiv0$).

А в вопросе о положительно определенных функциях указанного вида получен критерий, т.е.
необходимое и достаточное условие одновременно. Ответ полностью сводится к проверке
положительной определенности некоторой функции $f_{1}$ на прямой. Приведены примеры.

Первый вопрос --- это вопрос о принадлежности $f$ винеровской алгебре $\Big(x=(x_{1},
x_{2}, ..., x_{d}),$ $\ (x,y)=\sum\limits_{k=1}^{d}x_{k}y_{k},$ $\
|x|=\sqrt{(x,x)}\Big)$

\begin{equation*}
    A\big(\mathbb{R}^{d}\big)=\Big\{f(x)=\int\limits_{\mathbb{R}^{d}}g(y)e^{-i(x,y)}dy, \quad \|f\|_{A}=\int\limits_{\mathbb{R}^{d}}|g(y)|dy<\infty\Big\}.
\end{equation*}

Ответ на этот вопрос важен, напр., в проблеме мультипликаторов из
$L_{1}(\mathbb{R}^{d})$ в $L_{1}(\mathbb{R}^{d})$ (см. [\ref{Stein_Weiss}],
\textbf{3.19}).

Второй вопрос --- это вопрос о принадлежности алгебре
\begin{equation*}
A^{*}\big(\mathbb{R}^{d}\big)=\Big\{f(x)=\int\limits_{\mathbb{R}^{d}}g(y)e^{-i(x,y)}dy,
\quad \|f\|_{A^{*}}=\int\limits_{0}^{\infty}t^{d-1}\underset{|y|\geq
t}{esssup}|g(y)|dt<\infty\Big\}.
\end{equation*}

Эта алгебра при $d=1$ появилась в статье [\ref{Beurling}], а ее свойства и применение
см. в [\ref{Bel_Lifl_Trigub}].

Кроме того, воспользуемся обозначением
\begin{equation*}
A^{+}\big(\mathbb{R}^{d}\big)=\Big\{f(x)=\int\limits_{\mathbb{R}^{d}}g(y)e^{-i(x,y)}dy,
\quad g(y)\geq0,\quad y\in\mathbb{R}^{d}\Big\}.
\end{equation*}

См. также обзорную статью [\ref{Lifl_Samko_Trigub}].

Первые два вопроса тесно связаны, напр., с проблемами суммирования рядов и интегралов
Фурье по норме $L_{1}$ (или $C$) и почти всюду (точнее, в точках Лебега),
соответственно.

Точка Лебега $x_{0}$ для локально интегрируемой функции определяется следующим образом:
существует число $S(x_{0})$ такое, что
\begin{equation*}
    \lim\limits_{r\rightarrow0}\frac{1}{r^{d}}\int\limits_{|y|\leq
    r}\Big|f(x_{0}+y)-S(x_{0})\Big|dy=0.
\end{equation*}

А в третьем вопросе наша цель: свести положительную определенность функций
рассматриваемого вида к одномерному случаю. О подобных вопросах для радиальных функций,
т.е. функций, зависящих только от $|x|$, см. замечания 3 и 4 ниже.

Теперь о суммируемости рядов Фурье.

Если $W$ --- выпуклое и симметричное относительно нуля тело в $\mathbb{R}^{d}$, то
частные суммы ряда Фурье функций $f\in L_{1}(\mathbb{T}^{d})$, $\mathbb{T}=[-\pi,\pi]$,
$d\geq2$.

\begin{equation}\label{eq_1}
   S_{n}(f,W)=\sum\limits_{\frac{k}{n}\in W}\widehat{f}_{k}e_{k},\quad
   e_{k}=e^{i(k,x)},\quad\widehat{f}_{k}=\frac{1}{(2\pi)^{d}}\int\limits_{\mathbb{T}^{d}}f(x)e^{-i(k,x)}dx
\end{equation}
имеют, как известно, наибольший рост норм
\begin{equation*}
    \|S_{n}(\cdot,W)\|_{L_{1}\rightarrow L_{1}}=\frac{1}{(2\pi)^{d}}\int\limits_{\mathbb{T}^{d}}\Bigg|\sum\limits_{\frac{k}{n}\in W}e^{i(k,x)}\Bigg|dx
\end{equation*}
при $n\rightarrow\infty$ в случае шара $W\big(n^{\frac{d-1}{2}}\big)$, а наименьший ---
в случае куба $\big(\ln^{d}n\big)$. А с этим ростом связаны и разные требования
гладкости функций для того чтобы в $L_{1}(\mathbb{T}^{d})$
\begin{equation*}
    \lim\limits_{n\rightarrow\infty}S_{n}(f,W)=f
\end{equation*}
(см. [\ref{Trigub_Bel}], \textbf{9.2}).

Фейер доказал при $d=1$, что средние арифметические первых $n$ частных сумм сходятся на
всем пространстве $L_{1}(\mathbb{T})$, а Лебег доказал сходимость этих средних во всех
точках Лебега (почти всюду).

При $d\geq2$ Бохнер изучил сходимость средних (Бохнера -- Рисса), связанных с шаровыми
суммами. Для этих средних имеет место и сходимость во всех точках Лебега (см., напр.,
[\ref{Trigub_Bel}], с.~354).

Впервые Марцинкевич [\ref{Marcinkiewicz}] рассмотрел средние арифметические квадратных
частных сумм ($d=2$) и доказал их сходимость на $L_{1}(\mathbb{T}^{2})$. Как стало ясно
позже, этот результат эквивалентен тому, что
$\big(1-\max\big\{|x_{1}|,|x_{2}|\big\}\big)_{+}\in A(\mathbb{R}^{2})$. Отметим, что
сходимость таких средних на $L_{1}(\mathbb{T}^{2})$ эквивалентна сходимости на всем
пространстве $C(\mathbb{T}^{2})$ непрерывных и периодических функций. И это общий факт
для свёрточных операторов.

Изучение сходимости разных средних при $d=2$ почти всюду (без указания множества
сходимости) продолжено в работах [\ref{Gigiashv}], [\ref{Skopina}]. По поводу сходимости
на $L_{1}(\mathbb{T}^{d})$, $d\geq2$, см. [\ref{Podkor}] -- [\ref{Belinsky}].

В настоящей статье доказаны следующие теоремы, первые две из которых имеют общий
характер.

\begin{Th}\label{Th_1}
    Если $f\in L_{1}(\mathbb{R}^{d})$, то для того чтобы $f\in A(\mathbb{R}^{d})$
    необходимо и достаточно, чтобы $f\in C(\mathbb{R}^{d})$ и $\widehat{f}\in
    L_{1}(\mathbb{R}^{d})$.
\end{Th}

\begin{Th}\label{Th_2}
    Если тригонометрический ряд $\sum\limits_{k\in \mathbb{Z}^{d}}\varphi(k)e_{k}$ $\
    \big(e_{k}=e^{i(k,x)}\big)$ является рядом Фурье функции $\Phi\in L_{1}(\mathbb{T}^{d}),\
    \mathbb{T}=[-\pi,\pi]$, то
    \begin{equation*}
        \frac{1}{(2\pi)^{d}}\int\limits_{\mathbb{T}^{d}}\big|\Phi(x)\big|dx=\min\limits_{\varphi_{c}}\|\varphi_{c}\|_{A(\mathbb{R}^{d})},
    \end{equation*}
где $\varphi_{c}(k)=\varphi(k)$ для всех $k\in \mathbb{Z}^{d}$. И такое продолжение
$\varphi$ с минимальной нормой до целой функции экспоненциального типа не выше $\pi$ по
$x_{1}, x_{2}, ..., x_{d}$ существует единственное.
\end{Th}

Это еще одна связь норм свёрточных операторов в $L_{1}(\mathbb{T}^{d})$ (и
$C(\mathbb{T}^{d})$) с преобразованием Фурье.

\begin{Th}{\rm (достаточное условие принадлежности $A(\mathbb{R}^{2})$)}\label{Th_3}

    Пусть $f(x_{1},x_{2})=f_{0}\big(\max\{|x_{1}|,|x_{2}|\}\big)\in C\cap
    L_{1}(\mathbb{R}^{2})$, а $f_{0}\in AC_{loc}(0,+\infty)$. Если еще
    \begin{equation*}
        \int\limits_{0}^{1}\big|f'_{0}(t)\big|\ln^{2}\frac{2}{t}dt+\int\limits_{1}^{\infty}t^{2}\big|f'_{0}(t)\big|dt<\infty,
    \end{equation*}
    \begin{equation*}
        \int\limits_{0}^{1}\frac{\omega\big(f'_{0};t\big)_{1}}{t}\ln\frac{2}{t}dt+\int\limits_{0}^{1}\frac{\omega\big(f'_{1};t\big)_{1}}{t}dt<\infty,
    \end{equation*}
где $f'_{1}(t)=tf'_{0}(t)$, а $\omega(g;t)_{1}$ --- модуль непрерывности $g$ в метрике
$L_{1}(\mathbb{R}_{+})$, то $f\in A(\mathbb{R}^{2})$.
\end{Th}

Несколько более слабое достаточное условие применено к изучению разных методов
суммирования двойных рядов и интегралов Фурье в [\ref{Kotova_Trigub_UMG2}].

\begin{Th}{\rm (случай финитной и выпуклой функции $f_{0}$)}\label{Th_4}

    Если $f_{0}\in C(\mathbb{R}_{+})$, $f_{0}(t)=0$ при $t\geq1$, а на $[0,1]$ $f_{0}$
    убывает и выпукла вверх, то $f_{0}\big(\max\{|x_{1}|,|x_{2}|\}\big)\in
    A(\mathbb{R}^{2})$ тогда и только тогда, когда
\begin{equation*}
        \int\limits_{0}^{1}\frac{f_{0}(1-t)}{t}\ln\frac{2}{t}dt<\infty.
    \end{equation*}
\end{Th}

\begin{Th}\label{Th_5}
    Если $f_{0}\big(\max\{|x_{1}|,|x_{2}|\}\big)\in L_{1}(\mathbb{R}^{2})\cap A^{*}(\mathbb{R}^{2})$, $f_{0}\in
AC_{loc}(\mathbb{R}_{+})$
    и при некотором $\delta>0$
\begin{equation*}
        \int\limits_{0}^{\infty}e^{\delta t}\big|f'_{0}(t)\big|dt<\infty,
    \end{equation*}
то $f_{0}(t)\equiv0$.
\end{Th}

Из этой теоремы следует, что существует функция $f\in L_{1}(\mathbb{T}^{2})$, для
которой средние Марцинкевича, напр., не сходятся во всех ее точках Лебега
$\Big(\big(1-\max\{|x_{1}|,|x_{2}|\}\big)_{+}\notin A^{*}(\mathbb{R}^{2})\Big)$. А почти
всюду всегда сходятся [\ref{Gigiashv}]. Дело в том, что принадлежность $A^{*}$ --- это
главное достаточное условие для сходимости во всех точках Лебега ([\ref{Stein_Weiss}],
\textbf{1.25}). Но, как оказалось, это и необходимое условие (см. [\ref{Trigub_Bel}],
\textbf{8.1.3}).

\begin{Th}\label{Th_6}
    Пусть $f(x_{1},x_{2})=f_{0}\big(\max\{|x_{1}|,|x_{2}|\}\big)$, где $f_{0}\in C(\mathbb{R}_{+})\cap
AC_{loc}(0,+\infty)$ и $\lim\limits_{t\rightarrow+\infty}tf_{0}(t)=0$, а
$\int\limits_{0}^{\infty}t\big|f_{0}(t)\big|dt<\infty$. Если для всех
$y=(y_{1},y_{2})\in\mathbb{R}^{2}$
\begin{equation*}
    \widehat{f}(y)=\int\limits_{\mathbb{R}^{2}}f(x)e^{-i(x,y)}dx\geq0,
\end{equation*}
то и $\widehat{f}_{1}(x)\geq0$ для всех $x\in \mathbb{R}^{1}$, где
\begin{equation*}
    f_{1}(x)=|x|f_{0}\big(|x|\big)+\int\limits_{|x|}^{\infty}f_{0}(t)dt\qquad \Big(=-\int\limits_{|x|}^{\infty}tf'_{0}(t)dt\Big).
\end{equation*}
\end{Th}

\begin{Th}{\rm (в обозначениях теоремы \ref{Th_6})}\label{Th_7}

Если четная функция $f_{1}\in L_{1}(\mathbb{R}_{+})\cap C_{0}(\mathbb{R}_{+})\cap
AC_{loc}(0,+\infty)$,
\begin{equation*}
\int\limits_{0}^{\infty}(t+1)\big|f'_{1}(t)\big|dt<\infty
\end{equation*}
и сходится интеграл $\int\limits_{\rightarrow0}^{1}\frac{f'_{1}(t)}{t}dt$ (не
обязательно абсолютно), то из того, что $\widehat{f}_{1}(y)\geq0$ для всех $y\in
\mathbb{R}^{1}$, следует, что $\widehat{f}(y)\geq0$ для всех $y\in \mathbb{R}^{2}$.

А если $\widehat{f}_{1}(y)>0$ для всех $y\in \mathbb{R}^{1}$ и
$\int\limits_{0}^{\infty}tf'_{1}(t)dt<0$ (это и необходимо), то и $\widehat{f}(y)>0$ для
всех $y\in \mathbb{R}^{2}$.\\
\end{Th}

Добавление о строгой положительности преобразования Фурье сделано для того, чтобы можно
было применить аппроксимационную теорему Винера: для того чтобы линейные комбинации
сдвигов одной функции были плотны в $L_{1}$ необходимо и достаточно, чтобы ее
преобразование Фурье нигде не равнялось нулю при вещественных значениях переменных. План
статьи.

\S1 Доказательство теорем \ref{Th_1} и \ref{Th_2}.

\S2 Вспомогательные леммы.

\S3 Принадлежность $A(\mathbb{R}^{2})$ (теоремы \ref{Th_3} и \ref{Th_4}).

\S4 Принадлежность $A^{*}(\mathbb{R}^{2})$ (теорема \ref{Th_5}).

\S5 $A^{+}(\mathbb{R}^{2})$ (теоремы \ref{Th_6} и \ref{Th_7}).

\S6 Примеры. Связь с радиальными функциями.

\begin{center}
    \textbf{\S1 Доказательство теорем \ref{Th_1} и \ref{Th_2}}
\end{center}

Достаточность в теореме \ref{Th_1} сразу следует из формулы обращения преобразования
Фурье ([\ref{Stein_Weiss}], \textbf{1.21}).

Необходимость.

Если $g\in L_{1}(\mathbb{R}^{d})$, а преобразование Фурье ("нормированное")
\begin{equation*}
    \widehat{g}(y)=\frac{1}{(2\pi)^{\frac{d}{2}}}\int\limits_{\mathbb{R}^{d}}g(x)e^{-i(x,y)}dx,
\end{equation*}
то $\check{g}(y)=\widehat{g}(-y)$ --- обратное преобразование Фурье. Если еще $\widehat{g}\in
L_{1}(\mathbb{R}^{d})$, то в силу формулы обращения почти всюду
\begin{equation*}
    g=\check{\widehat{g}}=\widehat{\check{g}}.
\end{equation*}

Докажем, что если
\begin{equation}\label{eq_2}
    f=\widehat{f}_{1},\quad f_{1}\in L_{1},\quad f_{1}=\check{f}_{2},\quad f_{2}\in L_{1},
\end{equation}
то $f=f_{2}$ почти всюду.

В силу формулы умножения ([\ref{Stein_Weiss}], \textbf{1.15}), примененной дважды в
предположении, что $g$ и $\widehat{g}\in L_{1}$,
\begin{equation*}
    \int\limits_{\mathbb{R}^{d}}fg=\int\limits_{\mathbb{R}^{d}}\widehat{f}_{1}g=\int\limits_{\mathbb{R}^{d}}f_{1}\widehat{g}=\int\limits_{\mathbb{R}^{d}}\check{f}_{2}\widehat{g}=\int\limits_{\mathbb{R}^{d}}f_{2}g.
\end{equation*}

Так что для всех таких функций $g$
\begin{equation*}
    \int\limits_{\mathbb{R}^{d}}(f-f_{2})g=0,
\end{equation*}
где $f\in C_{0}(\mathbb{R}^{d})$, а $f_{2}\in L_{1}(\mathbb{R}^{d})$.

Докажем, что $f=f_{2}$ почти всюду на любом кубе $[a,b]^{d}$ ($b-a>2$).

А для этого достаточно показать, что
\begin{equation*}
\int\limits_{[a,b]^{d}}(f-f_{2})=0,
\end{equation*}
так как тогда при $d=2$, напр.,
\begin{equation*}
\int\limits_{a}^{t}dx_{1}\int\limits_{a}^{t}\big(f(x_{1},x_{2})-f_{2}(x_{1},x_{2})\big)dx_{2}=0
\end{equation*}
и после дифференцирования по $t$ дважды $f=f_{2}$ почти всюду.

Если при $n\in \mathbb{N}$ $\ h_{n}\in C(\mathbb{R})$, $h_{n}=1$ на
$\Big[a+\frac{1}{n},b-\frac{1}{n}\Big]$, $h_{n}=0$ на $\mathbb{R}\setminus[a,b]$ и
линейная на $\Big[a,a+\frac{1}{n}\Big]$ и $\Big[b-\frac{1}{n},b\Big]$, то в качестве
аппроксимации функции $g_{0}$, которая равна единице на $[a,b]^{2}$ и нулю на
$\mathbb{R}^{2}\setminus[a,b]^{2}$, возьмем последовательность функций
$g_{n}(x_{1},x_{2})=h_{n}(x_{1})\cdot h_{n}(x_{2})$. Очевидно, что $g_{n}$ и
$\widehat{g}_{n}\in L_{1}(\mathbb{R}^{2})$. Поэтому
\begin{equation}\label{eq_3}
\int\limits_{[a,b]^{2}}(f-f_{2})g_{n}=\int\limits_{\mathbb{R}^{2}}(f-f_{2})g_{n}=0.
\end{equation}

Кроме того,
\begin{equation*}
\Big|\int\limits_{[a,b]^{2}}f(g_{0}-g_{n})\Big|\leq\|f\|_{C[a,b]^{2}}\int\limits_{[a,b]^{2}}\big|g_{0}-g_{n}\big|=O\Big(\frac{1}{n}\Big),
\end{equation*}
и
\begin{equation*}
\Big|\int\limits_{[a,b]^{2}}f_{2}(g_{0}-g_{n})\Big|\leq\int\limits_{[a,b]^{2}}|f_{2}|\cdot\big|g_{0}-g_{n}\big|.
\end{equation*}

По теореме Лебега о мажорируемой сходимости, учитывая, что $\big|g_{0}-g_{n}\big|\leq1$,
последний интеграл стремится к нулю при $n\rightarrow\infty$. После перехода в
\eqref{eq_3} к пределу, получаем, что $\int\limits_{[a,b]^{2}}(f-f_{2})=0$ и, значит,
$f=f_{2}$ почти всюду.

Применяем доказанное утверждение \eqref{eq_2} к функции $F=\widehat{f}$:
\begin{equation*}
F=\widehat{f},\quad f\in L_{1},\quad f=\widehat{g},\quad g\in L_{1}\qquad \Rightarrow\qquad
F=\widehat{f}=g.
\end{equation*}

Следовательно, $\widehat{f}\in L_{1}$. Теорема \ref{Th_1} доказана.\\

\textbf{Доказательство теоремы \ref{Th_2}}. Поскольку коэффициенты Фурье равны
($k\in\mathbb{Z}^{d}$)
\begin{equation*}
    \varphi(k)=\frac{1}{(2\pi)^{d}}\int\limits_{\mathbb{T}^{d}}\Phi(x)e^{-i(k,x)}dx,
\end{equation*}
то полагая при $y\in \mathbb{R}^{d}$
\begin{equation*}
    \varphi_{c_{0}}(y)=\frac{1}{(2\pi)^{d}}\int\limits_{\mathbb{T}^{d}}\Phi(x)e^{-i(y,x)}dx,
\end{equation*}
получаем, что при любом продолжении функции $\varphi$ с $\mathbb{Z}^{d}$ имеем
\begin{equation*}
    \inf\limits_{\varphi_{c}}\|\varphi_{c}\|_{A}\leq\|\varphi_{c_{0}}\|_{A}=\frac{1}{(2\pi)^{d}}\int\limits_{\mathbb{T}^{d}}|\Phi(x)|dx.
\end{equation*}

С другой стороны, если
\begin{equation*}
    \varphi_{c}(y)=\int\limits_{\mathbb{R}^{d}}g(x)e^{-i(x,y)}dx,\quad
    \|\varphi_{c}\|_{A}=\|g\|_{L_{1}(\mathbb{R}^{d})}<\infty,
\end{equation*}
то при $k\in \mathbb{Z}^{d}$
\begin{equation*}
\begin{split}
    \varphi(k)=\varphi_{c}(k)&=\int\limits_{\mathbb{R}^{d}}g(x)e^{-i(k,x)}dx=\sum\limits_{m\in\mathbb{Z}^{d}}\int\limits_{\mathbb{T}^{d}}g(x)e^{-i(k,x-2\pi
    m)}dx=\\
    &=\sum\limits_{m\in\mathbb{Z}^{d}}\int\limits_{\mathbb{T}^{d}}g(x+2\pi
    m)e^{-i(k,x)}dx=\int\limits_{\mathbb{T}^{d}}g_{T}(x)e^{-i(k,x)}dx,
\end{split}
\end{equation*}
где $g_{T}(x)=\sum\limits_{m\in\mathbb{Z}^{d}}g(x+2\pi m)$,
$\|g_{T}\|_{L_{1}(\mathbb{T}^{d})}\leq\|g\|_{L_{1}(\mathbb{R}^{d})}$.

Здесь применена теорема Беппо Леви. Ряд для $g_{T}$ сходится почти всюду абсолютно.

Поскольку
\begin{equation*}
    g_{T}\sim\frac{1}{(2\pi)^{d}}\sum\limits_{k\in\mathbb{Z}^{d}}\varphi(k)e_{k}
\end{equation*}
(ряд Фурье), то при любом продолжении $\varphi_{c}$ имеем
\begin{equation*}
    \|\varphi_{c_{0}}\|_{A}=\frac{1}{(2\pi)^{d}}\int\limits_{k\in\mathbb{T}^{d}}|\Phi(x)|dx=\|g_{T}\|_{L_{1}(\mathbb{T}^{d})}\leq\|g\|_{L_{1}(\mathbb{R}^{d})}=\|\varphi_{c}\|_{A}.
\end{equation*}

С учетом нижеследующего замечания теорема \ref{Th_2} доказана.

\emph{\textbf{Замечание 1.}} \emph{Таких продолжений $\varphi_{c}$, имеющих минимальную
норму в $A(\mathbb{R}^{d})$, может быть много. Но существует единственное такое
продолжение до целой функции экспоненциального типа не выше $\pi$ по $x_{1}, x_{2},
...,x_{d}$ (см. выше $\varphi_{c_{0}}$).}

Нужно доказать, что если $\varphi$ такая функция и $\varphi(k)=0$ для всех $k\in
\mathbb{Z}^{d}$, то $\varphi\equiv0$.

При $d=1$ рассмотрим функцию $\psi(z)=\frac{\varphi(z)}{\sin\pi z}$ и убедимся в ее
ограниченности на всей плоскости $\mathbb{C}$.

По условию $|\varphi(z)|\leq Me^{\pi|Im z|}$ ($z\in\mathbb{C}$), так как $\varphi\in
B_{\pi}$.

А если $\pi|Im z|\geq\frac{1}{2}$, то при $Im z<0$, напр.,
\begin{equation*}
    |\sin\pi z|=\frac{1}{2}\big|e^{i\pi z}-e^{-i\pi z}\big|\geq\frac{1}{2}|e^{i\pi z}|\cdot\big(1-|e^{-2i\pi
    z}|\big)\geq\frac{1}{2}\Big(1-\frac{1}{e}\Big)\big|e^{i\pi
    z}\big|>\frac{1}{3}e^{\pi|Im z|}.
\end{equation*}

И функция $|\psi(z)|\leq3M$ при $\pi|Im z|\geq\frac{1}{2}$.

А в полосе $\pi|Im z|\leq\frac{1}{2}$ в силу периодичности $|\sin\pi z|$ можно
ограничиться одним прямоугольником при $|Re z|\leq\frac{1}{2}$ и применить принцип
максимума модуля.

Но ограниченная целая функция $\psi$ по теореме Лиувилля есть константа, т.е.
$\varphi(x)=c\sin\pi x$. А по лемме Римана--Лебега $\varphi(\infty)=0$. Так что и
$\varphi\equiv0$.

В кратном случае делим $\varphi$ на произведение синусов и применяем теорему Лиувилля и
лемму Римана--Лебега (см., напр., [\ref{Trigub_Bel}], \textbf{3.4.9}). Отметим еще, что,
как видно из доказательства, теорема \ref{Th_2} верна и для рядов Фурье любой конечной
борелевской комплекснозначной меры на торе. Вместо $A(\mathbb{R}^{d})$ --- алгебра
$B(\mathbb{R}^{d})$ преобразований Фурье мер:
\begin{equation*}
    B(\mathbb{R}^{d})=\Big\{f(x)=\int\limits_{\mathbb{R}^{d}}e^{-i(x,y)}d\mu(y),\quad \|f\|_{B}=var
    \mu<\infty\Big\}.
\end{equation*}

\emph{\textbf{Следствие.}} \emph{Если $\varphi\in B^{+}(\mathbb{R}^{d})$
($A^{+}(\mathbb{R}^{d})$), то для любого $\delta\neq0$
$\sum\limits_{k\in\mathbb{Z}^{d}}\varphi(\delta k)e_{k}$ --- ряд Фурье положительной на
торе $\mathbb{T}^{d}$ меры (функции).}

\textbf{Доказательство} при $d=1$, напр.,
\begin{equation*}
    \begin{split}
       &  \frac{1}{2\pi}var\mu=\frac{1}{2\pi}\sup\limits_{n}\int\limits_{-\pi}^{\pi}\Big|\sum\limits_{k}\varphi(\delta k)\Big(1-\frac{|k|}{n}\Big)_{+}e^{ikt}\Big|dt\leq\Big\|\Big(1-\frac{(\cdot)}{n}\Big)_{+}\varphi(\delta\cdot)\Big\|_{A}\leq\\
       &\leq\Big\|\Big(1-\frac{(\cdot)}{n}\Big)_{+}\Big\|_{A}\cdot\big\|\varphi(\delta\cdot)\big\|_{B^{+}}=1\cdot\varphi(0)=\frac{1}{2\pi}\int\limits_{-\pi}^{\pi}d\mu=\frac{1}{2\pi}\mu.
    \end{split}
\end{equation*}

\begin{center}
    \textbf{\S2 Вспомогательные леммы}
\end{center}

Предполагаем, что
\begin{equation}\label{eq_4}
f\in C(\mathbb{R}^{2})\cap L_{1}(\mathbb{R}^{2})\quad\textrm { и }\quad
f(x_{1},x_{2})=f_{0}\big(\max\big\{|x_{1}|,|x_{2}|\big\}\big).
\end{equation}

\begin{Lemma}\label{Lemma_1}
   При любом $p>0$
   \begin{equation*}
    \underset{\mathbb{R}^{2}}{\int\int}\Big|f(x_{1},x_{2})\Big|^{p}dx_{1}dx_{2}=8\int\limits_{0}^{\infty}t\big|f_{0}(t)\big|^{p}dt.
   \end{equation*}

Если $\widehat{f}_{0}(t)=\int\limits_{0}^{\infty}f_{0}(u)\sin utdu$ (синус--преобразование
Фурье), то при $y_{1}$ и $y_{2}\neq0$
\begin{equation*}
 \begin{split}
    \frac{1}{2}\widehat{f}(y_{1},y_{2})&=\frac{1}{2}\int\limits_{\mathbb{R}^{2}}f(x)e^{-i(x,y)}dx=\Big(\frac{1}{y_{1}}+\frac{1}{y_{2}}\Big)\widehat{f}_{0}(y_{1}+y_{2})-\Big(\frac{1}{y_{1}}-\frac{1}{y_{2}}\Big)\widehat{f}_{0}(y_{2}-y_{1})=\\
    &=-\frac{2}{y_{1}y_{2}}\int\limits_{\rightarrow0}^{\rightarrow\infty}f'_{0}(t)\sin
    ty_{1}\sin ty_{2}dt
  \end{split}
\end{equation*}
(последнее равенство при дополнительном предположении, что $f_{0}\in
AC_{loc}(0,+\infty)$).
\end{Lemma}

\textbf{Доказательство.} В силу симметрии $f$ относительно осей координат и биссектрисы
первой четверти
\begin{equation*}
\int\limits_{\mathbb{R}^{2}}\big|f(x)\big|^{p}dx=4\int\limits_{0}^{\infty}\int\limits_{0}^{\infty}\big|f(x_{1},x_{2})\big|^{p}dx_{1}dx_{2}=8\int\limits_{0}^{\infty}dx_{1}\int\limits_{0}^{x_{1}}\big|f_{0}(x_{1})\big|^{p}dx_{2}=8\int\limits_{0}^{\infty}x_{1}\big|f_{0}(x_{1})\big|^{p}dx_{1}.
\end{equation*}

Из того, что $f\in L_{1}(\mathbb{R}^{2})$ и доказанного равенства при $p=1$ следует, что
$tf_{0}(t)\in L_{1}(\mathbb{R}_{+})$. А так как $f\in C(\mathbb{R}^{2})$, то и $f_{0}\in
C(\mathbb{R}_{+})\cap L_{1}(\mathbb{R}_{+})$. Аналогично предыдущему
\begin{equation*}
    \begin{split}
       & \int\limits_{\mathbb{R}^{2}}f(x)e^{-i(x,y)}dx=4\int\limits_{0}^{\infty}\int\limits_{0}^{\infty}f_{0}\big(\max\big\{x_{1},x_{2}\big\}\big)\cos x_{1}y_{1}\cos x_{2}y_{2}dx_{1}dx_{2}= \\
       & =4\int\limits_{0}^{\infty}\cos x_{1}y_{1}dx_{1}\int\limits_{0}^{x_{1}}f_{0}(x_{1})\cos
       x_{2}y_{2}dx_{2}+4\int\limits_{0}^{\infty}\cos x_{2}y_{2}dx_{2}\int\limits_{0}^{x_{2}}f_{0}(x_{2})\cos
       x_{1}y_{1}dx_{1}=\\
       &=4\int\limits_{0}^{\infty}f_{0}(x_{1})\cos x_{1}y_{1}\frac{\sin x_{1}y_{2}}{y_{2}}dx_{1}+4\int\limits_{0}^{\infty}f_{0}(x_{2})\cos x_{2}y_{2}\frac{\sin x_{2}y_{1}}{y_{1}}dx_{2}=\\
       &=\frac{2}{y_{2}}\int\limits_{0}^{\infty}f_{0}(x_{1})\Big(\sin x_{1}(y_{1}+y_{2})+\sin
       x_{1}(y_{2}-y_{1})\Big)dx_{1}+\\
       &\quad+\frac{2}{y_{1}}\int\limits_{0}^{\infty}f_{0}(x_{2})\Big(\sin x_{2}(y_{1}+y_{2})-\sin
       x_{2}(y_{2}-y_{1})\Big)dx_{2}=\\
       &=2\Big(\frac{1}{y_{1}}+\frac{1}{y_{2}}\Big)\int\limits_{0}^{\infty}f_{0}(t)\sin
       t (y_{1}+y_{2})dt+2\Big(\frac{1}{y_{2}}-\frac{1}{y_{1}}\Big)\int\limits_{0}^{\infty}f_{0}(t)\sin
       t (y_{2}-y_{1})dt.
    \end{split}
\end{equation*}

Первое равенство доказано.

Интегрируем по частям, учитывая, что $f_{0}(+\infty)=0$ в силу леммы Римана--Лебега.
\begin{equation*}
    \begin{split}
       &\int\limits_{\mathbb{R}^{2}}f(x)e^{-i(x,y)}dx=2\Big(\frac{1}{y_{1}}+\frac{1}{y_{2}}\Big)f_{0}(t)\frac{-\cos
       t(y_{1}+y_{2})}{y_{1}+y_{2}}\Bigg|_{0}^{\infty}+2\Big(\frac{1}{y_{2}}-\frac{1}{y_{1}}\Big)\frac{-\cos
       t(y_{2}-y_{1})}{y_{2}-y_{1}}\Bigg|_{0}^{\infty}+\\
       &+2\Big(\frac{1}{y_{1}}+\frac{1}{y_{2}}\Big)\frac{1}{y_{1}+y_{2}}\int\limits_{0}^{\infty}f'_{0}(t)\cos
       t (y_{1}+y_{2})dt+2\Big(\frac{1}{y_{2}}-\frac{1}{y_{1}}\Big)\frac{1}{y_{2}-y_{1}}\int\limits_{0}^{\infty}f'_{0}(t)\cos
       t (y_{2}-y_{1})dt=\\
       &=\frac{2}{y_{1}y_{2}}\int\limits_{0}^{\infty}f'_{0}(t)\Big(\cos
       (y_{1}+y_{2})t-\cos
       (y_{2}-y_{1})t\Big)dt.
    \end{split}
\end{equation*}

Осталось применить известную формулу из тригонометрии.
\begin{Lemma}\label{Lemma_2}
   Если дополнительно к \eqref{eq_4} $\lim\limits_{t\rightarrow+\infty}tf_{0}(t)=0$,
   $f_{0}\in AC_{loc}(0,+\infty)$ и
   $\int\limits_{0}^{\infty}t\big|f'_{0}(t)\big|dt<\infty$, то при
   \begin{equation*}
    f_{1}(t)=tf_{0}(t)+\int\limits_{t}^{\infty}f_{0}(u)du=-\int\limits_{t}^{\infty}uf'_{0}(u)du\quad
    (t>0)
   \end{equation*}
   \begin{equation*}
    t\int\limits_{0}^{\infty}f_{1}(u)\cos utdu=\Big(t\widehat{f}_{0}(t)\Big)'.
   \end{equation*}
\end{Lemma}

\textbf{Доказательство.} Имеем, последовательно интегрируя по частям,
\begin{equation*}
    \begin{split}
       & \Big(t\widehat{f}_{0}(t)\Big)'=\widehat{f}_{0}(t)+t\widehat{f}'_{0}(t)=\int\limits_{0}^{\infty}f_{0}(u)\sin utdu+t\int\limits_{0}^{\infty}uf_{0}(u)\cos utdu= \\
       &=\int\limits_{0}^{\infty}f_{0}(u)\sin utdu+\Big[uf_{0}(u)\sin
       ut\Big]_{0}^{\infty}-\int\limits_{0}^{\infty}\Big(f_{0}(u)+uf'_{0}(u)\Big)\sin
       utdu=\\
       &=-\int\limits_{0}^{\infty}uf'_{0}(u)\sin
       utdu=\Big[\int\limits_{u}^{\infty}vf'_{0}(v)dv\cdot\sin
       ut\Big]_{0}^{\infty}-t\int\limits_{0}^{\infty}\cos ut
       du\int\limits_{u}^{\infty}vf'_{0}(v)dv=\\
       &=t\int\limits_{0}^{\infty}f_{1}(u)\cos ut du.
    \end{split}
\end{equation*}

При втором интегрировании по частям внеинтегральный член пропал, т.к.
\begin{equation*}
    u\int\limits_{u}^{\infty}vf'_{0}(v)dv=u\Big[vf_{0}(v)\Big]_{u}^{\infty}-\int\limits_{u}^{\infty}f_{0}(v)dv=-u^{2}f_{0}(u)-u\int\limits_{u}^{\infty}f_{0}(v)dv
\end{equation*}
и при $u\rightarrow+0$ стремится к нулю.

\begin{Lemma}\label{Lemma_3}
   В предположениях леммы \ref{Lemma_2} при некоторой абсолютной константе $C$
   \begin{equation*}
    \|\widehat{f}\|_{L_{1}(\mathbb{R}^{2})}\leq C
    \int\limits_{0}^{\infty}\Big|\Big(t\widehat{f}_{0}(t)\Big)'\Big|dt.
   \end{equation*}

\end{Lemma}

\textbf{Доказательство.} Пусть $g(t)=t\widehat{f}_{0}(t)$. В силу леммы \ref{Lemma_1}
\begin{equation*}
    \widehat{f}(y_{1},y_{2})=\frac{1}{y_{1}y_{2}}\Big(g(y_{1}+y_{2})-g(y_{2}-y_{1})\Big).
\end{equation*}

Используя симметрию $f$, получаем
\begin{equation*}
    \begin{split}
       & \|\widehat{f}\|_{L_{1}(\mathbb{R}^{2})}=8\int\limits_{0}^{\infty}y_{2}\int\limits_{0}^{y_{2}}\big|\widehat{f}(y_{1},y_{2})\big|dy_{1}\leq \\
       &\leq8\int\limits_{0}^{\infty}\frac{dy_{2}}{y_{2}}\int\limits_{0}^{y_{2}}\frac{1}{y_{1}}\Big|\int\limits_{y_{2}-y_{1}}^{y_{2}+y_{1}}g'(t)dt\Big|dy_{1}\leq8\int\limits_{0}^{\infty}\frac{dy_{2}}{y_{2}}\int\limits_{0}^{y_{2}}\frac{dy_{1}}{y_{1}}\int\limits_{y_{2}-y_{1}}^{y_{2}+y_{1}}|g'(t)|dt.
    \end{split}
\end{equation*}

Изменим порядок интегрирования:
\begin{equation*}
\begin{split}
&8\int\limits_{0}^{\infty}\frac{dy_{2}}{y_{2}}\int\limits_{0}^{2y_{2}}|g'(t)|dt\int\limits_{|y_{2}-t|}^{y_{2}}\frac{dy_{1}}{y_{1}}=8\int\limits_{0}^{\infty}\frac{dy_{2}}{y_{2}}\int\limits_{0}^{2y_{2}}|g'(t)|\ln\frac{y_{2}}{|y_{2}-t|}dt=\\
&=8\int\limits_{0}^{\infty}|g'(t)|dt\int\limits_{\frac{t}{2}}^{\infty}\frac{1}{y_{2}}\ln\frac{y_{2}}{|y_{2}-t|}dy_{2}.
\end{split}
\end{equation*}

Но внутренний интеграл ограничен по $t$ и $y_{2}$.

Для доказательства этого представим его в виде суммы трех интегралов по
$\Big[\frac{t}{2},t\Big]$, $[t,2t]$ и $[2t,+\infty)$. Первый интеграл не больше
\begin{equation*}
\frac{2}{t}\int\limits_{\frac{t}{2}}^{t}\ln\frac{y_{2}}{t-y_{2}}dy_{2}=\frac{2}{t}\Big[y_{2}\ln
y_{2}+(t-y_{2})\ln(t-y_{2})\Big]_{\frac{t}{2}}^{t}=2\ln2.
\end{equation*}
Второй интеграл не больше
\begin{equation*}
\frac{1}{t}\int\limits_{t}^{2t}\ln\frac{y_{2}}{y_{2}-t}dy_{2}=\frac{1}{t}\Big[y_{2}\ln
y_{2}-(y_{2}-t)\ln(y_{2}-t)\Big]_{t}^{2t}=2\ln2.
\end{equation*}
И третий интеграл не больше ($\ln(1+x)\leq x$ при $x\geq0$)
\begin{equation*}
\int\limits_{2t}^{\infty}\frac{1}{y_{2}}\ln\Big(1+\frac{t}{y_{2}-t}\Big)dy_{2}\leq\int\limits_{2t}^{\infty}\Big(\frac{1}{y_{2}-t}-\frac{1}{y_{2}}\Big)dy_{2}=\Big[\ln\frac{y_{2}-t}{y_{2}}\Big]_{2t}^{\infty}=\ln2
\end{equation*}

Лемма доказана.

Заметим, что если $f_{0}(t)=(1-t)_{+}$, то
$\int\limits_{0}^{\infty}\Big|\Big(t\widehat{f}_{0}(t)\Big)'\Big|dt=\infty$.

\begin{Lemma}{\rm (необходимое условие для принадлежности $A(\mathbb{R}^{2})$)}\label{Lemma_4}

   Если $f\in A(\mathbb{R}^{2})$ и $supp f\subset[-1,1]^{2}$, то сходится интеграл
   \begin{equation*}
      \int\limits_{\rightarrow0}^{1}\frac{f_{0}(1-t)}{t}\ln\frac{2}{t}dt.
   \end{equation*}

\end{Lemma}

\textbf{Доказательство.} Начнем с общего необходимого условия (для всех функций из
$A(\mathbb{R}^{2})$). Пусть смешанная разность равна
\begin{equation*}
    \triangle_{h_{1},h_{2}}^{1,1}\
    f(x_{1},x_{2})=f(x_{1}+h_{1},x_{2}+h_{2})-f(x_{1}+h_{1},x_{2}-h_{2})-f(x_{1}-h_{1},x_{2}+h_{2})+f(x_{1}-h_{1},x_{2}-h_{2}).
\end{equation*}
Тогда
\begin{equation*}
    \triangle_{h_{1},h_{2}}^{1,1}\
    e^{i(x_{1}y_{1}+x_{2}y_{2})}=-4e^{i(x_{1}y_{1}+x_{2}y_{2})}\sin h_{1}y_{1}\sin
    h_{2}y_{2}.
\end{equation*}

Пусть $f=\widehat{g}$, где $g\in L_{1}(\mathbb{R}^{2})$. Тогда при $0<a_{k}<b_{k}$ ($k=1$ и
$2$)
\begin{equation*}
    \begin{split}
       & \int\limits_{a_{1}}^{b_{1}}du_{1}\int\limits_{a_{2}}^{b_{2}}\frac{1}{u_{1}u_{2}} \triangle_{h_{1},h_{2}}^{1,1}\ f(x_{1},x_{2})du_{2}=\\
       &= \int\limits_{a_{1}}^{b_{1}}\frac{du_{1}}{u_{1}}\int\limits_{a_{2}}^{b_{2}}\frac{du_{2}}{u_{2}}\int\limits_{-\infty}^{\infty}\int\limits_{-\infty}^{\infty}g(y_{1},y_{2}) \triangle_{u_{1},u_{2}}^{1,1}\
       e^{-i(x_{1}y_{1}+x_{2}y_{2})}dy_{1}dy_{2}=\\
       &=-4\int\limits_{a_{1}}^{b_{1}}\frac{du_{1}}{u_{1}}\int\limits_{a_{2}}^{b_{2}}\frac{du_{2}}{u_{2}}\int\limits_{-\infty}^{\infty}\int\limits_{-\infty}^{\infty}g(y_{1},y_{2}) e^{-i(x_{1}y_{1}+x_{2}y_{2})}\sin u_{1}y_{1}\sin
       u_{2}y_{2}dy_{1}dy_{2}=\\
       &=-4\int\limits_{-\infty}^{\infty}\int\limits_{-\infty}^{\infty}g(y_{1},y_{2}) e^{-i(x_{1}y_{1}+x_{2}y_{2})}dy_{1}dy_{2}\int\limits_{a_{1}}^{b_{1}}\int\limits_{a_{2}}^{b_{2}}\frac{\sin u_{1}y_{1}}{u_{1}}\cdot\frac{\sin
       u_{2}y_{2}}{u_{2}}du_{1}du_{2}.
    \end{split}
\end{equation*}

Но
\begin{equation*}
    \sup\limits_{a,b,y}\Big|\int\limits_{a}^{b}\frac{\sin uy}{u}du\Big|=\sup\Big|\int\limits_{ay}^{by}\frac{\sin
    t}{t}dt\Big|<\infty.
\end{equation*}

Поскольку $g\in L_{1}(\mathbb{R}^{2})$, то можно перейти к пределу при $a_{1}$ и
$a_{2}\rightarrow+0$, а $b_{1}$ и $b_{2}\rightarrow+\infty$. Получаем
\begin{equation*}
    \int\limits_{\rightarrow0}^{\rightarrow\infty}\int\limits_{\rightarrow0}^{\rightarrow\infty}\frac{\triangle_{h_{1},h_{2}}^{1,1}\
    f(x_{1},x_{2})}{u_{1}u_{2}}du_{1}du_{2}=-4f(x_{1},x_{2})\Big(\int\limits_{0}^{\infty}\frac{\sin u}{u}du\Big)^{2}.
\end{equation*}

Так что для любой функции $f\in A(\mathbb{R}^{2})$ в любой точке $(x_{1},x_{2})\in
\mathbb{R}^{2}$ сходится указанный интеграл (не обязательно абсолютно).

А теперь возьмем функцию вида $f_{0}\big(\max\big\{|x_{1}|,|x_{2}|\big\}\big)$ с
условием $f_{0}(t)=0$ при $t\geq1$ в точке $(x_{1},x_{2})=(1,1)$.

Сходится интеграл
\begin{equation*}
    \int\limits_{\rightarrow0}^{2}\int\limits_{\rightarrow0}^{2}g_{1}(u_{1},u_{2})du_{1}du_{2},\qquad
    g_{1}(u_{1},u_{2})=\frac{f_{0}\big(\max\big\{|1-u_{1}|,|1-u_{2}|\big\}\big)}{u_{1}u_{2}}.
\end{equation*}

Разобьем квадрат $[0,2]^{2}$ диагоналями на четыре треугольника и воспользуемся
симметрией $g_{1}$. Интеграл равен
\begin{equation*}
    2\int\limits_{0}^{1}du_{2}\int\limits_{u_{2}}^{2-u_{2}}g_{1}(u_{1},u_{2})du_{1}+2\int\limits_{1}^{2}du_{1}\int\limits_{2-u_{1}}^{u_{1}}g_{1}(u_{1},u_{2})du_{2}.
\end{equation*}

В первом интеграле $0\leq u_{2}\leq1$ и $u_{2}\leq u_{1}\leq2-u_{2}$, откуда следует,
что $|1-u_{1}|\leq|1-u_{2}|=1-u_{2}$. Поэтому первое слагаемое равно
\begin{equation*}
2\int\limits_{0}^{1}du_{2}\int\limits_{u_{2}}^{2-u_{2}}\frac{f_{0}(1-u_{2})}{u_{1}u_{2}}du_{1}=2\int\limits_{0}^{1}\frac{f_{0}(1-u_{2})}{u_{2}}\ln\Big(\frac{2-u_{2}}{u_{2}}\Big)du_{2}.
\end{equation*}

Учтем теперь, что при $u\in (0,1]$
\begin{equation*}
    \Big|\ln\Big(\frac{2-u}{u}\Big)-ln\frac{2}{u}\Big|=\Big|\ln\Big(1-\frac{u}{2}\Big)\Big|\leq u\ln2.
\end{equation*}

Модуль второго слагаемого после грубой оценки не больше
\begin{equation*}
    \begin{split}
       & 2\sup\limits_{[0,1]}|f_{0}(t)|\int\limits_{1}^{2}du_{1}\int\limits_{2-u_{1}}^{u_{1}}\frac{du_{2}}{u_{1}u_{2}}=2\sup\limits_{[0,1]}|f_{0}(t)|\int\limits_{1}^{2}\frac{1}{u_{1}}\ln\frac{u_{1}}{2-u_{1}}du_{1}\leq \\
       &
       \leq2\sup\limits_{[0,1]}|f_{0}(t)|\int\limits_{1}^{2}\ln\frac{u_{1}}{2-u_{1}}du_{1}=2\ln2\cdot\sup\limits_{[0,1]}|f_{0}(t)|.
    \end{split}
\end{equation*}

Лемма доказана.

\begin{Lemma}\label{Lemma_5}

При $\alpha\in(-1,0]$ и $x\rightarrow+\infty$
\begin{equation*}
\int\limits_{0}^{1}(1-t)^{\alpha}e^{itx}dt=\frac{\Gamma(\alpha+1)}{x^{1+\alpha}}e^{i\big(x-\frac{\pi(\alpha+1)}{2}\big)}+\frac{i}{x}-\frac{\alpha}{x^{2}}+O\Big(\frac{|\alpha|}{x^{3}}\Big).
\end{equation*}

Для определения асимптотики при $\alpha>0$ нужно воспользоваться формулой ($x>0$):
\begin{equation*}
\int\limits_{0}^{1}(1-t)^{\alpha}e^{itx}dt=\frac{i}{x}-\frac{\alpha
i}{x}\int\limits_{0}^{1}(1-t)^{\alpha-1}e^{itx}dt.
\end{equation*}

\end{Lemma}

\textbf{Доказательство.} При $\alpha\in(-1,0)$ воспользуемся известным интегралом
\begin{equation*}
\int\limits_{0}^{\infty}t^{\alpha}e^{-it}dt=\Gamma(\alpha+1)e^{-i\frac{\pi(\alpha+1)}{2}}.
\end{equation*}

Тогда при $x>0$ (замена $x(1-t)=u$)
\begin{equation*}
\int\limits_{0}^{1}(1-t)^{\alpha}e^{itx}dt=\frac{e^{ix}}{x^{1+\alpha}}\int\limits_{0}^{x}u^{\alpha}e^{-iu}du=\frac{e^{ix}}{x^{1+\alpha}}\Big(\Gamma(\alpha+1)e^{-i\frac{\pi(\alpha+1)}{2}}-\int\limits_{x}^{\infty}u^{\alpha}e^{-iu}du\Big).
\end{equation*}

Осталось учесть, что
\begin{equation*}
    \begin{split}
     & \int\limits_{x}^{\infty}u^{\alpha}e^{-iu}du=\Big[ie^{-iu}u^{\alpha}\Big]_{x}^{\infty}-\alpha
    i\int\limits_{x}^{\infty}u^{\alpha-1}e^{-iu}du=-ie^{-ix}x^{\alpha}-\alpha i\Big(\Big[ie^{-iu}u^{\alpha-1}\Big]_{x}^{\infty}- \\
     &-(\alpha-1)i\int\limits_{x}^{\infty}u^{\alpha-2}e^{-iu}du\Big)=-ie^{-ix}x^{\alpha}-\alpha
     e^{-ix}x^{\alpha-1}-\alpha(\alpha-1)\int\limits_{x}^{\infty}u^{\alpha-2}e^{-iu}du=\\
     &=-ie^{-ix}x^{\alpha}-\alpha
     e^{-ix}x^{\alpha-1}-\alpha(\alpha-1)\Big(-ie^{-ix}x^{\alpha-2}-(\alpha-2)i\int\limits_{x}^{\infty}u^{\alpha-3}e^{-iu}du\Big)=\\
     &=-ie^{-ix}x^{\alpha}-\alpha
     e^{-ix}x^{\alpha-1}+i\alpha(\alpha-1)e^{-ix}x^{\alpha-2}-i\alpha(\alpha-1)(\alpha-2)\int\limits_{x}^{\infty}u^{\alpha-3}e^{-iu}du.
\end{split}
\end{equation*}

Можно выписать весь асимптотический ряд при $x\rightarrow+\infty$. Получаем
\begin{equation*}
\int\limits_{0}^{1}(1-t)^{\alpha}e^{itx}dt=\frac{\Gamma(\alpha+1)}{x^{1+\alpha}}e^{i\big(x-\frac{\pi(\alpha+1)}{2}\big)}+\frac{i}{x}-\frac{\alpha}{x^{2}}+O\Big(\frac{|\alpha|}{x^{3}}\Big).
\end{equation*}

\emph{\textbf{Замечание 2.}} \emph{Вместо интегрирования по частям иногда можно
воспользоваться разложением экспоненты в степенной ряд и свойствами $B$--функции
Эйлера.}

Например, при $\alpha>-1$ и $x>0$
\begin{equation*}
    (\alpha+1)x^{\alpha+1}\int\limits_{0}^{1}(1-t)^{\alpha}\sin txdt=\Big(x^{\alpha+2}\int\limits_{0}^{1}(1-t)^{\alpha+1}\sin txdt\Big)'.
\end{equation*}

\begin{center}
    \textbf{\S3 Принадлежность $A(\mathbb{R}^{2})$ (теоремы \ref{Th_3} и \ref{Th_4})}
\end{center}

Достаточное условие содержится в теореме \ref{Th_3} (см. во введении).

\textbf{Доказательство теоремы \ref{Th_3}.} Нужно проверить, что $\widehat{f}\in
L_{1}(\mathbb{R}^{2})$ (см. теорему \ref{Th_1}).

Применяем лемму \ref{Lemma_1} и неравенство $\big|\sin t\big|\leq|t|$. При
$(y_{1},y_{2})\in \mathbb{R}^{2}$
\begin{equation*}
    \big|\widehat{f}(y_{1},y_{2})\big|\leq4\int\limits_{0}^{\infty}t^{2}\big|f'_{0}(t)\big|dt\leq4\int\limits_{0}^{1}\big|f'_{0}(t)\big|dt+4\int\limits_{1}^{\infty}t^{2}\big|f'_{0}(t)\big|dt<\infty.
\end{equation*}

Поэтому интеграл по квадрату $[-1,1]^{2}$ конечен. Оставшуюся в силу симметрии область,
определяемую неравенствами $y_{1}>0$ и $y_{2}\geq\max\{1,y_{1}\}$, разобьем на три:
\begin{equation}\label{eq_5}
    D_{1}:\ y_{1}\in(0,1],\ y_{2}\geq1;\quad D_{2}:\ 1\leq y_{1}\leq y_{2}\leq
    y_{1}+1;\quad D_{3}:\ 1\leq y_{1}\leq y_{1}+1\leq
    y_{2}.
\end{equation}

Отметим сначала, что при $h\in L_{1}(\mathbb{R}_{+})$ и $t>0$
\begin{equation}\label{eq_6}
    \Big|\int\limits_{0}^{\infty}h(u)e^{iut}du\Big|\leq\frac{1}{2}\int\limits_{0}^{\infty}\Big|h(u)-h\Big(u+\frac{\pi}{t}\Big)\Big|du+\frac{1}{2}\int\limits_{0}^{\frac{\pi}{t}}|h(u)|du.
\end{equation}

Действительно,
\begin{equation*}
\int\limits_{0}^{\infty}h(u)e^{iut}du=-\int\limits_{-\frac{\pi}{t}}^{\infty}h\Big(u+\frac{\pi}{t}\Big)e^{iut}du=-\int\limits_{0}^{\infty}h\Big(u+\frac{\pi}{t}\Big)e^{iut}du+\int\limits_{0}^{\frac{\pi}{t}}h(u)e^{iut}du,
\end{equation*}
откуда следует, что
\begin{equation*}
\int\limits_{0}^{\infty}h(u)e^{iut}du=\frac{1}{2}\int\limits_{0}^{\infty}\Big(h(u)-h\Big(u+\frac{\pi}{t}\Big)\Big)e^{iut}du+\frac{1}{2}\int\limits_{0}^{\frac{\pi}{t}}h(u)e^{iut}du.
\end{equation*}

Осталось учесть, что $\big|e^{iut}\big|=1$.

\textbf{Случай $D_{1}$.} В силу леммы \ref{Lemma_1} и формулы интегрирования по частям
\begin{equation*}
\begin{split}
    &\widehat{f}(y_{1},y_{2})=-\frac{4}{y_{1}y_{2}}\int\limits_{0}^{\infty}f'_{0}(t)\sin
    ty_{1}\sin ty_{2}dt=\\
    &=-\frac{4}{y_{1}y_{2}}\Big(-\Big[\sin ty_{1}\int\limits_{t}^{\infty}f'_{0}(u)\sin
    uy_{2}du\Big]_{0}^{\infty}+\int\limits_{0}^{\infty}y_{1}\cos y_{1}tdt\int\limits_{t}^{\infty}f'_{0}(u)\sin
    uy_{2}du\Big).
    \end{split}
\end{equation*}

Внеинтегральный член равен нулю в силу того, что $f'_{0}\in L_{1}(\mathbb{R}_{+})$.
\begin{equation*}
\begin{split}
&\int\limits_{t}^{\infty}f'_{0}(u)\sin
uy_{2}du=\frac{1}{2}\int\limits_{t}^{\infty}f'_{0}(u)e^{iuy_{2}}du-\frac{1}{2}\int\limits_{t}^{\infty}f'_{0}(u)e^{-iuy_{2}}du=\\
&=\frac{1}{2}e^{ity_{2}}\int\limits_{0}^{\infty}f'_{0}(u+t)e^{iuy_{2}}du-\frac{1}{2}e^{-ity_{2}}\int\limits_{0}^{\infty}f'_{0}(u+t)e^{-iuy_{2}}du.
\end{split}
\end{equation*}

Применяем неравенство \eqref{eq_6}:
\begin{equation*}
\Big|\int\limits_{t}^{\infty}f'_{0}(u)\sin
uy_{2}du\Big|\leq\frac{1}{2}\int\limits_{0}^{\infty}\Big|f'_{0}(u+t)-f'_{0}\Big(u+t+\frac{\pi}{y_{2}}\Big)\Big|du+\frac{1}{2}\int\limits_{0}^{\frac{\pi}{y_{2}}}\Big|f'_{0}(u+t)\Big|du.
\end{equation*}

Следовательно,
\begin{equation*}
\big|\widehat{f}(y_{1},y_{2})\big|\leq\frac{2}{y_{2}}\int\limits_{0}^{\infty}dt\int\limits_{t}^{\infty}\Big|f'_{0}(u)-f'_{0}\Big(u+\frac{\pi}{y_{2}}\Big)\Big|du+\frac{2}{y_{2}}\int\limits_{0}^{\infty}dt\int\limits_{t}^{t+\frac{\pi}{y_{2}}}\Big|f'_{0}(u)\Big|du.
\end{equation*}

Второе слагаемое после перемены порядка интегрирования равно
\begin{equation*}
\frac{2}{y_{2}}\int\limits_{0}^{\frac{\pi}{y_{2}}}u\big|f'_{0}(u)\big|du+\frac{2\pi}{y_{2}^{2}}\int\limits_{\frac{\pi}{y_{2}}}^{\infty}\big|f'_{0}(u)\big|du\leq\frac{2\pi}{y_{2}^{2}}\int\limits_{0}^{\infty}\big|f'_{0}(u)\big|du.
\end{equation*}

А первое слагаемое не больше
\begin{equation*}
    \begin{split}
       & \frac{2}{y_{2}}\int\limits_{0}^{\infty}\Big|uf'_{0}(u)-\Big(u+\frac{\pi}{y_{2}}\Big)f'_{0}\Big(u+\frac{\pi}{y_{2}}\Big)\Big|du+\frac{2\pi}{y_{2}^{2}}\int\limits_{0}^{\infty}\big|f'_{0}(u)\big|du= \\
       &=\frac{2}{y_{2}}\int\limits_{0}^{\infty}\Big|f'_{1}(u)-f'_{1}\Big(u+\frac{\pi}{y_{2}}\Big)\Big|du+O\Big(\frac{1}{y^{2}_{2}}\Big).
    \end{split}
\end{equation*}

Поэтому
\begin{equation*}
    \underset{D_{1}}{\int\int}\big|\widehat{f}(y_{1},y_{2})\big|dy_{1}dy_{2}\leq\int\limits_{0}^{1}dy_{1}\int\limits_{1}^{\infty}\frac{2}{y_{2}}\omega\Big(f'_{1};\frac{\pi}{y_{2}}\Big)_{1}dy_{2}+\int\limits_{1}^{\infty}O\Big(\frac{1}{y_{2}^{2}}\Big)dy_{2}.
\end{equation*}

Но
\begin{equation*}
\int\limits_{1}^{\infty}\frac{\omega\Big(f'_{1};\frac{\pi}{y_{2}}\Big)_{1}}{y_{2}}dy_{2}=\int\limits_{0}^{\pi}\frac{\omega(f'_{1};t)_{1}}{t}dt<\infty.
\end{equation*}

Так что
\begin{equation*}
\underset{D_{1}}{\int\int}\big|\widehat{f}(y_{1},y_{2})\big|dy_{1}dy_{2}<\infty.
\end{equation*}

\textbf{Случай $D_{2}$} (см. \eqref{eq_5}). При $t>0$
\begin{equation*}
    \widehat{f}_{0}(t)=\int\limits_{0}^{\infty}f_{0}(u)\sin
    utdu=\frac{f_{0}(0)}{t}+\frac{1}{t}\int\limits_{0}^{\infty}f'_{0}(u)\cos
    utdu.
\end{equation*}

К интегралу применяем неравенство \eqref{eq_6}.
\begin{equation*}
    \begin{split}
       & \Big|t\widehat{f}_{0}(t)-f_{0}(0)\Big|\leq\frac{1}{2}\Big|\int\limits_{0}^{\infty}f'_{0}(u)e^{iut}du\Big|+\frac{1}{2}\Big|\int\limits_{0}^{\infty}f'_{0}(u)e^{-iut}du\Big|\leq \\
       &\leq\frac{1}{2}\int\limits_{0}^{\infty}\Big|f'_{0}(u)-f'_{0}\Big(u+\frac{\pi}{t}\Big)\Big|du+\frac{1}{2}\int\limits_{0}^{\frac{\pi}{t}}\big|f'_{0}(u)\big|du.
    \end{split}
\end{equation*}

В силу леммы \ref{Lemma_1}
\begin{equation*}
    \widehat{f}(y_{1},y_{2})=\frac{1}{y_{1}y_{2}}\Big((y_{1}+y_{2})\widehat{f}_{0}(y_{1}+y_{2})-(y_{2}-y_{1})\widehat{f}_{0}(y_{2}-y_{1})\Big).
\end{equation*}

Поэтому, учитывая, как и раньше, определение $\omega(g;h)_{1}$, имеем
\begin{equation*}
\begin{split}
&\underset{D_{2}}{\int\int}\big|\widehat{f}(y_{1},y_{2})\big|dy_{1}dy_{2}=\int\limits_{1}^{\infty}dy_{1}\int\limits_{y_{1}}^{y_{1}+1}\big|\widehat{f}(y_{1},y_{2})\big|dy_{2}\leq\\
&\leq\frac{1}{2}\int\limits_{1}^{\infty}\frac{dy_{1}}{y_{1}}\int\limits_{y_{1}}^{y_{1}+1}\frac{dy_{2}}{y_{2}}\Big(\omega\Big(f'_{0};\frac{\pi}{y_{1}+y_{2}}\Big)_{1}+\omega\Big(f'_{0};\frac{\pi}{y_{2}-y_{1}}\Big)_{1}+\\
&+\int\limits_{0}^{\frac{\pi}{y_{1}+y_{2}}}\big|f'_{0}(u)\big|du+\int\limits_{0}^{\frac{\pi}{y_{2}-y_{1}}}\big|f'_{0}(u)\big|du\Big)\leq\int\limits_{1}^{\infty}\frac{dy_{1}}{y_{1}}\int\limits_{y_{1}}^{y_{1}+1}\frac{\omega\Big(f'_{0};\frac{\pi}{y_{2}-y_{1}}\Big)_{1}}{y_{2}}dy_{2}+\\
&+\int\limits_{1}^{\infty}\frac{dy_{1}}{y_{1}}\int\limits_{y_{1}}^{y_{1}+1}\frac{dy_{2}}{y_{2}}\int\limits_{0}^{\frac{\pi}{y_{2}-y_{1}}}\big|f'_{0}(u)\big|du.
\end{split}
\end{equation*}

Первый интеграл равен $\big(\omega(g;h)_{1}\leq2\|g\|_{L_{1}}\big)$
\begin{equation*}
\int\limits_{1}^{\infty}\frac{dy_{1}}{y_{1}}\int\limits_{0}^{1}\frac{dy_{2}}{y_{1}+y_{2}}\omega\Big(f'_{0};\frac{\pi}{y_{2}}\Big)_{1}\leq2\int\limits_{0}^{\infty}\big|f'_{0}(t)\big|dt\int\limits_{1}^{\infty}\frac{dy_{1}}{y_{1}^{2}}=2\int\limits_{0}^{\infty}\big|f'_{0}(t)\big|dt.
\end{equation*}

И второй интеграл не больше
\begin{equation*}
\int\limits_{0}^{\infty}\big|f'_{0}(u)\big|du\int\limits_{1}^{\infty}\frac{dy_{1}}{y_{1}}\int\limits_{y_{1}}^{y_{1}+1}\frac{dy_{2}}{y_{2}}\leq\int\limits_{0}^{\infty}\big|f'_{0}(u)\big|du.
\end{equation*}

Так что
\begin{equation*}
\underset{D_{2}}{\int\int}\big|\widehat{f}(y_{1},y_{2})\big|dy_{1}dy_{2}<3\int\limits_{0}^{\infty}\big|f'_{0}(u)\big|du.
\end{equation*}

\textbf{Случай $D_{3}$} (см. \eqref{eq_5}). Как и в предыдущем случае
\begin{equation*}
\begin{split}
&\underset{D_{2}}{\int\int}\big|\widehat{f}(y_{1},y_{2})\big|dy_{1}dy_{2}=\int\limits_{1}^{\infty}dy_{1}\int\limits_{y_{1}+1}^{\infty}\big|\widehat{f}(y_{1},y_{2})\big|dy_{2}\leq\\
&\leq\int\limits_{1}^{\infty}\frac{dy_{1}}{y_{1}}\int\limits_{y_{1}+1}^{\infty}\frac{\omega\Big(f'_{0};\frac{\pi}{y_{2}-y_{1}}\Big)_{1}}{y_{2}}dy_{2}+\int\limits_{1}^{\infty}\frac{dy_{1}}{y_{1}}\int\limits_{y_{1}+1}^{\infty}\frac{dy_{2}}{y_{2}}\int\limits_{0}^{\frac{\pi}{y_{2}-y_{1}}}\big|f'_{0}(u)\big|du.
\end{split}
\end{equation*}

Первый интеграл равен
\begin{equation*}
\begin{split}
&\int\limits_{1}^{\infty}\frac{dy_{1}}{y_{1}}\int\limits_{1}^{\infty}\frac{\omega\Big(f'_{0};\frac{\pi}{y_{2}}\Big)_{1}}{y_{2}+y_{1}}dy_{2}=\int\limits_{1}^{\infty}\omega\Big(f'_{0};\frac{\pi}{y_{2}}\Big)_{1}dy_{2}\int\limits_{1}^{\infty}\frac{dy_{1}}{y_{1}(y_{1}+y_{2})}=\\
&=\int\limits_{1}^{\infty}\omega\Big(f'_{0};\frac{\pi}{y_{2}}\Big)_{1}\ln(1+y_{2})\frac{dy_{2}}{y_{2}}=\int\limits_{0}^{\pi}\frac{\omega(f'_{0};u)_{1}}{u}\ln\Big(1+\frac{\pi}{u}\Big)du<\infty.
\end{split}
\end{equation*}

Второй интеграл равен
\begin{equation*}
\begin{split}
&\int\limits_{1}^{\infty}\frac{dy_{1}}{y_{1}}\int\limits_{1}^{\infty}\frac{dy_{2}}{y_{1}+y_{2}}\int\limits_{0}^{\frac{\pi}{y_{2}}}\big|f'_{0}(u)\big|du=\int\limits_{1}^{\infty}dy_{2}\int\limits_{0}^{\frac{\pi}{y_{2}}}\big|f'_{0}(u)\big|du\int\limits_{1}^{\infty}\frac{dy_{1}}{y_{1}(y_{1}+y_{2})}=\\
&=\int\limits_{1}^{\infty}\frac{1}{y_{2}}\ln(1+y_{2})dy_{2}\int\limits_{0}^{\frac{\pi}{y_{2}}}\big|f'_{0}(u)\big|du=\int\limits_{0}^{\pi}\big|f'_{0}(u)\big|du\int\limits_{1}^{\frac{\pi}{u}}\frac{\ln(1+y_{2})}{y_{2}}dy_{2}\leq\\
&\leq2\int\limits_{0}^{\pi}\big|f'_{0}(u)\big|du\int\limits_{1}^{\frac{\pi}{u}}\frac{\ln(1+y_{2})}{1+y_{2}}dy_{2}=\int\limits_{0}^{\pi}\big|f'_{0}(u)\big|\ln^{2}\Big(1+\frac{\pi}{u}\Big)du<\infty.
\end{split}
\end{equation*}

Теорема \ref{Th_3} доказана.
\\

Рассмотрим теперь случай финитной и выпуклой функции $f_{0}$.

\textbf{Доказательство теоремы \ref{Th_4}.} Применяем теорему \ref{Th_3}. $f_{0}(t)=0$
при $t\geq1$. А так как правая производная $f'_{0}(0)\in[-f_{0}(0),0]$, то первое
интегральное условие теоремы \ref{Th_3} с $f'_{0}$ выполнено.

В силу монотонности $f'_{0}$ на $(0,1)$
\begin{equation*}
    \begin{split}
       & \omega\big(f'_{0};h\big)_{1}=\sup\limits_{0<\delta\leq h\leq1}\int\limits_{0}^{\infty}\big|f'_{0}(t)-f'_{0}(t+\delta)\big|dt\leq\sup\Big|\int\limits_{0}^{1-\delta}\big(f'_{0}(t)-f'_{0}(t+\delta)\big)dt\Big|+ \\
       &+\sup\Big|\int\limits_{1-\delta}^{1}f'_{0}(t)dt\Big|=\sup\big|f_{0}(1-\delta)-f_{0}(0)-f_{0}(1)+f_{0}(\delta)\big|+\sup\big|f_{0}(1)-f_{0}(1-\delta)\big|\leq\\
       &\leq\sup\big|f_{0}(0)-f_{0}(\delta)\big|+2\sup\big|f_{0}(1)-f_{0}(1-\delta)\big|.
    \end{split}
\end{equation*}

Функция $f_{0}$ убывает на $[0,1]$ и $f_{0}(1)=0$. Поэтому
\begin{equation*}
    \omega\big(f'_{0};h\big)_{1}\leq2f_{0}(1-h)+f_{0}(0)-f_{0}(h)\leq2f_{0}(1-h)+f'_{0}(0)h
\end{equation*}
и
\begin{equation*}
    \int\limits_{0}^{1}\frac{\omega\big(f'_{0};t\big)_{1}}{t}\ln\frac{2}{t}dt\leq2\int\limits_{0}^{1}\frac{f_{0}(1-t)}{t}\ln\frac{2}{t}dt+f'_{0}(0)\int\limits_{0}^{1}\ln\frac{2}{t}dt<\infty.
\end{equation*}

Проверим теперь условие на функцию $f'_{1}(t)=tf'_{0}(t)$.

При $0<\delta\leq h\leq1$
\begin{equation*}
    \begin{split}
       & \int\limits_{0}^{1}\big|tf'_{0}(t)-(t+\delta)f'_{0}(t+\delta)\big|dt\leq\delta\int\limits_{0}^{1}\big|f'_{0}(t+\delta)\big|dt+ \\
       &+\Big|\int\limits_{0}^{1-\delta}t\big(f'_{0}(t)-f'_{0}(t+\delta)\big)dt\Big|+\int\limits_{1-\delta}^{1}t\big|f'_{0}(t)\big|dt\leq\\
       &\leq\delta\int\limits_{0}^{\infty}\big|f'_{0}(t)\big|dt+\Bigg|\Big[t\big(f_{0}(t)-f_{0}(t+\delta)\big)\Big]_{0}^{1-\delta}-\int\limits_{0}^{1-\delta}\big(f_{0}(t)-f_{0}(t+\delta)\big)dt\Bigg|+\\
       &+\int\limits_{1-\delta}^{1}\big|f'_{0}(t)\big|dt\leq h
       \int\limits_{0}^{1}\big|f'_{0}(t)\big|dt+3f_{0}(1-h).
    \end{split}
\end{equation*}
($\big|f(t)-f(t+\delta)\big|$ с ростом $t$ и $\delta$ возрастает).

Следовательно, и
\begin{equation*}
\int\limits_{0}^{1}\frac{\omega\big(f'_{1};t)_{1}}{t}dt\leq\int\limits_{0}^{1}\big|f'_{0}(t)\big|dt+3\int\limits_{0}^{1}\frac{f_{0}(1-t)}{t}dt<\infty.
\end{equation*}

Достаточность доказана.

А необходимость содержится в лемме \ref{Lemma_4}.\\
\newpage
\emph{\textbf{Замечание 3.} }(о выпуклых функциях)

\emph{В двумерном радиальном случае
$f(x_{1},x_{2})=f_{0}\big(\sqrt{x_{1}^{2}+x_{2}^{2}}\big)$, где $f_{0}$ выпукла на
$[0,1]$, необходимое и достаточное условие принадлежности $A(\mathbb{R}^{2})$ выглядит
так:
\begin{equation*}
\int\limits_{0}^{1}\frac{f_{0}(1-t)}{t^{\frac{3}{2}}}dt<\infty.
\end{equation*}
(см. [\ref{Trigub_Bel}], \textbf{6.5.7}). В одномерном случае для преобразования Фурье
выпуклой на отрезке функции есть асимптотическая формула (см. там же \textbf{6.4.7}).
Впрочем, асимптотическая формула есть и для радиальных функций при $d\geq2$ с условием
выпуклости на производную (см. [\ref{Trigub_1977}]).}

А теорему \ref{Th_4} легко обобщить на кусочно выпуклые функции, используя локальное
свойство алгебры $A$.

\begin{center}
    \textbf{\S4 Принадлежность $A^{*}(\mathbb{R}^{2})$ (теорема \ref{Th_5})}
\end{center}

Для методов суммирования рядов Фурье, определяемых функцией--множителем с условием типа
выпуклости, вообще нет разницы между сходимостью на $L_{1}$ и во всех точках Лебега (см.
[\ref{Trigub_Bel}], c.~354--355). Другое дело в случае функции--множителе вида
$f_{0}\big(\max\big\{|x_{1}|,|x_{2}|\big\}\big)$. См. ниже пример.

\textbf{Доказательство теоремы \ref{Th_5}.} По условию $f\in L_{1}\cap
A^{*}(\mathbb{R}^{2})$. Пусть $f=\widehat{g}$. После перехода к полярным координатам
\begin{equation*}
    \|g\|_{L_{1}(\mathbb{R}^{2})}=\int\limits_{0}^{\infty}tdt\int\limits_{0}^{2\pi}\Big|g(t\cos\varphi,t\sin\varphi)\Big|d\varphi\leq2\pi\int\limits_{0}^{\infty}t\
    \underset{y_{1}^{2}+y_{2}^{2}\geq t^{2}}{esssup}\big|g(y_{1},y_{2})\big|dt<\infty.
\end{equation*}

По формуле обращения
\begin{equation*}
    f(x_{1},x_{2})=\frac{1}{4\pi^{4}}\underset{\mathbb{R}^{2}}{\int\int}\widehat{f}(y_{1},y_{2})e^{i(x_{1}y_{1}+x_{2}y_{2})}dy_{1}dy_{2}.
\end{equation*}

Равенство имеет место всюду в силу непрерывности $f$ и правой части.

В силу леммы \ref{Lemma_1}
\begin{equation*}
    \frac{1}{2}\widehat{f}(y_{1},y_{2})=\Big(\frac{1}{y_{1}}+\frac{1}{y_{2}}\Big)\widehat{f}_{0}(y_{1}+y_{2})-\Big(\frac{1}{y_{1}}-\frac{1}{y_{2}}\Big)\widehat{f}_{0}(y_{2}-y_{1}).
\end{equation*}

При любом $a>0$, $y_{1}=t>0$ и $y_{2}=a+t$
\begin{equation}\label{eq_7}
\begin{split}
    &\|f\|_{A^{*}}=\frac{1}{2\pi^{2}}\int\limits_{0}^{\infty}t\sup\limits_{y_{1}^{2}+y_{2}^{2}\geq
    t^{2}}\Big|\Big(\frac{1}{y_{1}}-\frac{1}{y_{2}}\Big)\widehat{f}_{0}(y_{2}-y_{1})-\Big(\frac{1}{y_{1}}+\frac{1}{y_{2}}\Big)\widehat{f}_{0}(y_{1}+y_{2})\Big|dt\geq\\
    &\geq\frac{1}{2\pi^{2}}\int\limits_{0}^{\infty}t\Big|\Big(\frac{1}{t}-\frac{1}{t+a}\Big)\widehat{f}_{0}(a)-\Big(\frac{1}{t}+\frac{1}{t+a}\Big)\widehat{f}_{0}(2t+a)\Big|dt\geq\\
    &\geq\frac{1}{2\pi^{2}}\int\limits_{1}^{\infty}\Big|\frac{a\widehat{f}_{0}(a)}{t+a}-\frac{2t+a}{t+a}\widehat{f}_{0}(2t+a)\Big|dt.
\end{split}
\end{equation}

Учтём теперь, что при $t>0$ и $t\rightarrow+\infty$
\begin{equation}\label{eq_8}
    \widehat{f}_{0}(t)=\frac{f_{0}(0)}{t}+\frac{1}{t}\int\limits_{0}^{\infty}f'_{0}(u)\cos
    utdu=\frac{f_{0}(0)}{t}+o\Big(\frac{1}{t}\Big).
\end{equation}

Предположим, что $\widehat{f}_{0}(a)\neq0$. Сравним сходимость (расходимость) интеграла
\eqref{eq_7} с более простым.
\begin{equation*}
    \lim\limits_{t\rightarrow+\infty}\Big(\frac{a\widehat{f}_{0}(a)}{t+a}-\frac{2t+a}{t+a}\widehat{f}_{0}(2t+a)\Big):\frac{a\widehat{f}_{0}(a)}{t+a}=1-\lim\limits_{t\rightarrow+\infty}\frac{(2t+a)\widehat{f}_{0}(2t+a)}{a\widehat{f}_{0}(a)}=1-\frac{f_{0}(0)}{a\widehat{f}_{0}(a)}.
\end{equation*}

Но
\begin{equation*}
\int\limits_{1}^{\infty}\frac{a\big|\widehat{f}_{0}(a)\big|}{t+a}dt=\infty,
\end{equation*}
а интеграл \eqref{eq_7} сходится. Следовательно, $a\widehat{f}_{0}(a)=f_{0}(0)$ и в силу
\eqref{eq_8}
\begin{equation*}
\int\limits_{0}^{\infty}f'_{0}(u)\cos audu=0.
\end{equation*}

Если этот интеграл равен нулю почти всюду или на плотном в $\mathbb{R}_{+}$ множестве,
то делаем вывод, что $\widehat{f}_{0}(a)=0$ на этом множестве, а в силу непрерывности ---
всюду и $f\equiv0$. Из дополнительного условия
$\int\limits_{0}^{\infty}\big|f'_{0}(u)\big|e^{\delta u}du<\infty$ следует, что
преобразование Фурье $\int\limits_{0}^{\infty}f'_{0}(u)\cos tudu$ является сужением на
вещественную ось аналитической в полосе $|Im\ z|<\delta$ функции. Поэтому нули этой
функции изолированные.

Теорема \ref{Th_5} доказана.\\

\textbf{\emph{Пример}.} Средние Марцинкевича--Рисса
\begin{equation*}
    \sum\limits_{k_{1}=-n}^{n}\sum\limits_{k_{2}=-n}^{n}\Big(1-\Big(\frac{1}{n}\max\big\{|k_{1}|,|k_{2}|\big\}\Big)^{\alpha}\Big)^{\beta}_{+}\widehat{f}_{k_{1},k_{2}}e^{i(k_{1}x_{1}+k_{2}x_{2})}
\end{equation*}
при любых $\alpha$ и $\beta>0$ сходятся при $n\rightarrow\infty$ на всем пространстве
$L_{1}(\mathbb{T}^{2})$ в силу теорем \ref{Th_2} и \ref{Th_3}. Но в силу теоремы
\ref{Th_5} существует функция $f\in L_{1}(\mathbb{R}^{2})$, для которой эти средние
расходятся хотя бы в одной точке Лебега.

\begin{center}
    \textbf{\S5 Положительно определенные функции ($A^{+}$)}
\end{center}

Проверка положительной определенности функции
$f_{0}\big(\max\big\{|x_{1}|,|x_{2}|\big\}\big)$ (принадлежность
$A^{+}(\mathbb{R}^{2})$) полностью сводится к проверке принадлежности
$A^{+}(\mathbb{R}^{1})$ функции
\begin{equation*}
    f_{1}(x)=|x|f_{0}\big(|x|\big)+\int\limits_{|x|}^{\infty}f_{0}(t)dt=-\int\limits_{|x|}^{\infty}tf'_{0}(t)dt.
\end{equation*}

\textbf{Доказательство теоремы \ref{Th_6}} (см. во введении).

В силу леммы \ref{Lemma_1} ($p=1$) $f\in L_{1}(\mathbb{R}^{2})$, а так как
$\widehat{f}(y_{1},y_{2})\geq0$ и $f$ ограничена около нуля, то и $\widehat{f}\in
L_{1}(\mathbb{R}^{2})$ (см. [\ref{Stein_Weiss}], \textbf{1.26} и [\ref{Trigub_Bel}],
\textbf{3.1.15}).

В силу условий теоремы и $f_{0}\in L_{1}(\mathbb{R}_{+})$.

По этой же лемме \ref{Lemma_1} при $y_{1}y_{2}\neq0$
\begin{equation}\label{eq_9}
    \frac{1}{2}\widehat{f}(y_{1},y_{2})=\frac{1}{y_{1}y_{2}}\Big((y_{1}+y_{2})\widehat{f}_{0}(y_{1}+y_{2})-(y_{2}-y_{1})\widehat{f}_{0}(y_{2}-y_{1})\Big),
\end{equation}
откуда видно, что для того чтобы $\widehat{f}(y_{1},y_{2})\geq0$ при $y_{2}\geq y_{1}>0$ не
только достаточно, но и необходимо, чтобы функция $t\widehat{f}_{0}(t)$ возрастала на
$\mathbb{R}_{+}$. Но если она возрастает на $\mathbb{R}_{+}$, то в силу нечетности
$\widehat{f}_{0}$ и $\widehat{f}(y_{1},y_{2})\geq0$ при $y_{1}>y_{2}>0$. Следовательно,
$\widehat{f}(y_{1},y_{2})\geq0$ при $y_{1}>0$ и $y_{2}>0$.

Но функция $\widehat{f}(y_{1},y_{2})$ четная по $y_{1}$ и $y_{2}$, поэтому
$\widehat{f}(y_{1},y_{2})\geq0$ почти всюду (кроме, может быть, двух осей координат). Но
$\widehat{f}\in C(\mathbb{R}^{2})$ и $\widehat{f}(y_{1},y_{2})\geq0$ на $\mathbb{R}^{2}$.

Учитывая еще, что $f_{0}(+\infty)=0$, получаем
\begin{equation*}
    g(t)=t\widehat{f}_{0}(t)=t\int\limits_{0}^{\infty}f_{0}(u)\sin utdu=f_{0}(0)+\int\limits_{0}^{\infty}f'_{0}(u)\cos
    utdu.
\end{equation*}

А эта дифференцируемая при $t>0$ функция возрастает в том и только в том случае, если ее
производная
\begin{equation*}
g'(t)=-\int\limits_{0}^{\infty}uf'_{0}(u)\sin utdu\geq0.
\end{equation*}

В силу леммы \ref{Lemma_2}
\begin{equation*}
g'(t)=t\int\limits_{0}^{\infty}f_{1}(u)\cos utdu.
\end{equation*}

Так что при $t>0$
\begin{equation*}
\int\limits_{0}^{\infty}f_{1}(u)\cos utdu\geq0.
\end{equation*}

При этом $f_{1}\in C(\mathbb{R}_{+})$ и
\begin{equation*}
\int\limits_{0}^{\infty}\big|f_{1}(t)\big|dt\leq\int\limits_{0}^{\infty}t\big|f_{0}(t)\big|dt+\int\limits_{0}^{\infty}dt\int\limits_{t}^{\infty}\big|f_{0}(u)\big|du=2\int\limits_{0}^{\infty}t\big|f_{0}(t)\big|dt<\infty,
\end{equation*}
т.е. $f_{1}\in A^{+}(\mathbb{R}^{1})$.

Теорема \ref{Th_6} доказана.\\

Это был переход от $f_{0}$ к $f_{1}$. А теперь перейдем в обратном направлении с
добавлением строгой положительности преобразование Фурье.

\textbf{Доказательство теоремы \ref{Th_7}.} Полагаем при $t>0$
\begin{equation*}
    f_{0}(t)=\int\limits_{t}^{\infty}\frac{f_{1}(u)}{u}du.
\end{equation*}
и проверим, что $f_{0}$ удовлетворяет условиям теоремы \ref{Th_6}.

Очевидно, что $f_{0}\in C(\mathbb{R}_{+})\cap C^{1}(0,+\infty)$.
\begin{equation*}
    \big|tf_{0}(t)\big|\leq
    t\int\limits_{t}^{\infty}\frac{1}{u}\big|f'_{1}(u)\big|du\leq\int\limits_{t}^{\infty}\big|f'_{1}(u)\big|du\rightarrow0\quad
    (t\rightarrow+\infty)
\end{equation*}
\begin{equation*}
\int\limits_{0}^{\infty}t\big|f_{0}(t)\big|dt\leq\int\limits_{0}^{\infty}tdt\int\limits_{t}^{\infty}\frac{1}{u}\big|f'_{1}(u)\big|du=\int\limits_{0}^{\infty}\big|f'_{1}(u)\big|\frac{du}{u}\int\limits_{0}^{u}tdt=\frac{1}{2}\int\limits_{0}^{\infty}u\big|f'_{1}(u)\big|du<\infty.
\end{equation*}

Кроме того,
\begin{equation*}
\int\limits_{0}^{\infty}t\big|f'_{0}(t)\big|dt=\int\limits_{0}^{\infty}\big|f'_{1}(t)\big|dt\leq\int\limits_{0}^{\infty}(t+1)\big|f'_{1}(t)\big|dt<\infty.
\end{equation*}

Теперь (см. доказательство теоремы \ref{Th_6}) из неравенства при $t>0$
\begin{equation*}
\int\limits_{0}^{\infty}f_{1}(u)\cos utdu\geq0\quad (>0)
\end{equation*}
следует, что при $t>0$
\begin{equation*}
-\int\limits_{0}^{\infty}uf'_{0}(u)\sin utdu\geq0\quad (>0)
\end{equation*}

Но тогда функция $g(t)=t\widehat{f}_{0}(t)$ возрастает (строго возрастает), так как
$g'(t)\geq0$ ($>0$) на $(0,+\infty)$.

Теперь из формулы \eqref{eq_9} выводим, что $\widehat{f}(y_{1},y_{2})\geq0$ ($>0$) при
$y_{1}y_{2}\neq0$. Из той же формулы \eqref{eq_9} следует, что при $y_{2}>0$
\begin{equation*}
    \widehat{f}(0,y_{2})=\lim\limits_{y_{1}\rightarrow0}\widehat{f}(y_{1},y_{2})=\frac{1}{y_{2}}\cdot\lim\limits_{y_{1}\rightarrow0}\frac{g(y_{1}+y_{2})-g(y_{2}-y_{1})}{y_{1}}=\frac{2}{y_{2}}g'(y_{2})>0.
\end{equation*}

Но $g'(0)=0$ и
\begin{equation*}
    \widehat{f}(0,0)=\lim\limits_{y_{2}\rightarrow0}\frac{2}{y_{2}}\big(g'(y_{2})-g'(0)\big)=2g''(0).
\end{equation*}

Осталось учесть, что $g'(t)=-\int\limits_{0}^{\infty}uf'_{0}(u)\sin utdu$ и
\begin{equation*}
    g''(0)=\lim\limits_{t\rightarrow0}\frac{g'(t)}{t}=-\int\limits_{0}^{\infty}u^{2}f'_{0}(u)du=-\int\limits_{0}^{\infty}uf'_{1}(u)du>0.
\end{equation*}

Теорема \ref{Th_7} доказана.

\begin{center}
    \textbf{\S6 Примеры. Связь с радиальными функциями}
\end{center}

Остановимся на применении теорем \ref{Th_6} и \ref{Th_7}. Дело в том, что примеров
положительно определенных радиальных функций накопилось уже довольно много.

\emph{\textbf{Пример 1.}} Проверим, что преобразование Фурье функции
$e^{-\max\big\{|x_{1}|,|x_{2}|\big\}}$ больше нуля всюду на $\mathbb{R}^{2}$.

В этом случае $f_{0}(t)=e^{-t}$, а $f_{1}(x)=\big(|x|+1\big)e^{-|x|}$
($x\in\mathbb{R}$).

А так как
\begin{equation*}
\int\limits_{0}^{\infty}(t+1)e^{-t}\cos
txdt=Re\int\limits_{0}^{\infty}(t+1)e^{t(ix-1)}dt=\frac{2}{(1+x^{2})^{2}}>0
\end{equation*}
и
\begin{equation*}
\int\limits_{0}^{\infty}tf'_{1}(t)dt=-\int\limits_{0}^{\infty}f_{1}(t)dt<0,
\end{equation*}
то можно применить теорему \ref{Th_7}.

Так что имеем аналог метода суммирования Абеля--Пуассона (свертка функции $f$ с
положительным ядром)
\begin{equation*}
    \sum\limits_{k_{1}=-\infty}^{\infty}\sum\limits_{k_{2}=-\infty}^{\infty}e^{-\varepsilon\max\big\{|k_{1}|,|k_{2}|\big\}}\widehat{f}_{k_{1},k_{2}}e^{i(k_{1}x_{1}+k_{2}x_{2})}\underset{\varepsilon\rightarrow+0}{\longrightarrow}f(x_{1},x_{2}).
\end{equation*}

Сходимость в $L_{1}(\mathbb{T}^{2})$ и в $C(\mathbb{T}^{2})$ имеет место всегда, но не
всегда в точках Лебега (см. теорему \ref{Th_5}).

Отметим теперь такой общий факт. Число производных при $x>0$ у $f_{0}$ и $f_{1}$
одинаковое, а из того, что $f_{0}\in C\cap L_{1}(\mathbb{R}_{+})$ следует, что
$f'_{1}(0)=0$.

Действительно, при $x\neq0$
\begin{equation*}
    \frac{f_{1}(x)-f_{1}(0)}{x}=f_{0}(x)-\frac{1}{x}\int\limits_{0}^{x}f_{0}(t)dt
\end{equation*}
и после перехода к пределу при $x\rightarrow0$ получаем $f'_{1}(0)=0$.

Теперь перейдем к радиальным финитным функциям.

Легко видеть, что для любой вещественной функции $g\in L_{1}\cap A(\mathbb{R}^{d})$ и
четной по $x_{1}, ..., x_{d}$ свертка
\begin{equation*}
    (g\ast g)(x)=\int\limits_{\mathbb{R}^{d}}g(y)g(x-y)dy\in L_{1}\cap A^{+}(\mathbb{R}^{d})
\end{equation*}
(преобразование Фурье свертки равно $(\widehat{g})^{2}$).

А поскольку преобразование Фурье радиальной функции --- радиальная функция, то и свертка
радиальная, если $g$ --- радиальная.

А если носитель $g$ в шаре $|x|\leq\sigma$, то носитель свертки в шаре $|x|\leq2\sigma$.
С таким представлением связан критерий М.~Г.~Крейна положительной определенности
финитных функций (см. [\ref{Lifl_Samko_Trigub}], \textbf{11.9}) и общий критерий
Ф.~Рисса принадлежности $A$ (свертка двух функций из $L_{2}$, [\ref{Lifl_Samko_Trigub}],
\textbf{4.10}).

Простейшими в некотором смысле являются функции вида $p_{n}(|x|)$ (алгебраический
полином степени $n$) при $|x|\leq1$ и $0$ при $|x|\geq1$. Еще в 1987 г. автор построил
такие функции $A_{r,d}\in C^{2r}(\mathbb{R}^{d})\cap A^{+}(\mathbb{R}^{d})$ при $r+1\in
\mathbb{N}$ и нечетной размерности $d$ наименьшей возможной степени
$n=3r+\frac{d+1}{2}$.

При условии $A_{r,d}(0)=1$ такие "сплайны" определяются однозначно (публикация в
[\ref{Trigub_Tezis}]). Дело в том, что мою статью с доказательством при $d\geq2$ тогда
отклонили в одном из московских журналов и пришлось публиковать в малодоступных изданиях
(см. ссылки в [\ref{Trigub_Bel}]). В 1995 г. к этим сплайнам другим путем пришел
H.~Wendland [\ref{Wendland}].

Используя однозначность определения таких сплайнов в $C^{2r}(\mathbb{R}^{d})$ и их
гладкость в нуле и на границе носителя, легко их найти. При $r\geq1$ и степени
$n=3r+\frac{d+1}{2}$
\begin{equation*}
    A_{r,d}(t)=(1-t)^{2r+\frac{d+1}{2}}_{+}\Big(1+\sum\limits_{k=1}^{r}a_{k}t^{k}\Big),
\end{equation*}
где коэффициенты $\{a_{k}\}_{1}^{r}$ определяются из условия, что у $A_{r,d}$ не должно
быть нечетных степеней $t=|x|$ до $2r-1$, т.е.
\begin{equation*}
    \frac{d^{2\nu+1}}{dt^{2\nu+1}}\big\{A_{r,d}(t)\big\}_{t=0}=0\quad (0\leq\nu\leq r-1)
\end{equation*}

Это линейная система уравнений с ненулевым определителем.

\emph{Примеры.}
\begin{equation}\label{eq_10}
\begin{split}
   &A_{1,d}(x)=\big(1-|x|\big)^{\frac{d+5}{2}}_{+}\Big(1+\frac{d+5}{2}|x|\Big),\\
   &A_{2,d}(x)=\big(1-|x|\big)^{\frac{d+7}{2}}_{+}\Big(1+\frac{d+7}{2}|x|+\frac{(d+5)(d+9)}{12}|x|^{2}\Big).
\end{split}
\end{equation}

При нечетном $d$ и $n\geq\frac{d+1}{2}$ существует сплайн указанного вида из
$C^{r}(\mathbb{R}^{d})$ при $r=2\Big[\frac{2n-d-1}{6}\Big]$ и только при указанном $r$
(не больше). Но при нечетном $r$ может не быть однозначности (см. теорему 1 в
[\ref{Trigub_2002}]). Обобщение $A_{r,d}$ на два нецелых параметра предложил
В.~П.~Заставный (см. [\ref{Zastavn_Trigub}]):
\begin{equation}\label{eq_11}
    h_{\mu,\nu}(x)=\big(1-|x|\big)^{\mu+\nu-1}_{+}\int\limits_{0}^{1}t^{\mu-1}(1-t)^{\nu-1}\big(1-t+(1+t)|x|\big)^{\nu-1}dt.
\end{equation}

При этом указаны все возможные значения параметров $\mu$ и $\nu$, при которых
$h_{\mu,\nu}\in A^{+}(\mathbb{R}^{1})$ и $\widehat{h}_{\mu,\nu}(x)>0$ на $\mathbb{R}^{1}$.
$h_{r+1,r+1}=cA_{r,1}$. А переход от $A^{+}(\mathbb{R}^{1})$ к $A^{+}(\mathbb{R}^{d})$
($d\geq2$) давно известен (см. ниже). А еще более общее (четырехпараметрическое)
семейство изучил M.~D.~Buhman [\ref{Buhmanm}] (см. также [\ref{Zastavn_2006}]).

А теперь вернемся к функциям вида
$f_{1}(x_{1},x_{2})=f_{0}\big(\max\big\{|x_{1}|,|x_{2}|\big\}\big)$, $f_{0}(t)=0$ при
$t\geq1$.

Начнем с необходимых условий.

При $tf'_{0}(t)\in AC[0,1]$ после двукратного интегрирования по часятм
\begin{equation*}
    \widehat{f}_{1}(x)=-2\int\limits_{0}^{1}\cos
    uxdu\int\limits_{u}^{1}tf'_{0}(t)dt=\frac{2}{x^{2}}\Big(f'_{0}(1)\cos
    x-\int\limits_{0}^{1}\big(tf'_{0}(t)\big)'\cos txdt\Big).
\end{equation*}

Если $f'_{0}(1)\neq0$, то $\widehat{f}_{1}$ меняет знак бесконечное число раз.

А если и $\Big(tf'_{0}(t)\Big)'\in AC[0,1]$, то и $f''_{0}(1)=0$, так как
\begin{equation*}
    \widehat{f}_{1}(x)=\frac{2}{x^{3}}\Big(-f''_{0}(1)\sin
    x+\int\limits_{0}^{1}\big(tf'_{0}(t)\big)'\sin txdt\Big).
\end{equation*}

\emph{\textbf{Пример 2.}} $f_{0}(t)=(1-t)^{\alpha}_{+}$, $\alpha\geq3$.

В силу указанных необходимых условий должно быть $\alpha>2$. А так как
\begin{equation*}
f_{1}(t)=c(1-t)^{\alpha}_{+}(1+\alpha t),
\end{equation*}
то при $\alpha\in(2,3)$ можно, используя лемму \ref{Lemma_5}, найти асимптотику
$\widehat{f}_{1}$:
\begin{equation*}
    \int\limits_{0}^{1}(1-t)^{\alpha}(1+\alpha t)\cos
    txdu=-\frac{\alpha^{2}(\alpha+1)}{x^{1+\alpha}}\Gamma(\alpha-2)\sin\Big(x-\pi\frac{\alpha-2}{2}\Big)+O\Big(\frac{1}{x^{4}}\Big).
\end{equation*}

Остается случай $\alpha\geq3$.

При $\alpha\in\mathbb{N}$ $\ f_{0}\big(\max\big\{|x_{1}|,|x_{2}|\big\}\big)\in
A^{+}(\mathbb{R}^{2})$ и преобразование Фурье больше нуля всюду, т.к.
\begin{equation*}
    \big(1-|x|\big)^{\alpha}_{+}\big(1+\alpha|x|\big)=A_{1,d}\big(|x|\big)\quad
    (d=2\alpha-5)
\end{equation*}
(см. \eqref{eq_10}), а $A^{+}(\mathbb{R}^{d})\subset A^{+}(\mathbb{R}^{1})$.

А при любом $\alpha\geq3$ $\ f_{1}(t)=\ch_{\alpha-1,2}(t)$ (см. \eqref{eq_11}).

\emph{\textbf{Пример 3.}} Если в качестве $f_{1}$ взять сплайн $A_{r,d}$ максимальной
гладкости (см. выше), то получим в квадрате $[-1,1]^{2}$ полином от
$\max\big\{|x_{1}|,|x_{2}|\big\}$, не содержащий максимальное число первых нечетных
степеней $\max\big\{|x_{1}|,|x_{2}|\big\}$.

Отметим еще, что если существует $f''_{1}(0)$, то при довольно общих условиях $f_{1}$ и
$(-f''_{1})$ являются положительно определенными функциями одновременно (см.
[\ref{Trigub_Bel}], \textbf{6.2.15}).

\textbf{\emph{Замечание 4.}} (о радиальных функциях)

\emph{Через $f_{d}:\mathbb{R}_{+}\rightarrow\mathbb{R}$ обозначим такую функцию, что
$f_{d}\big(|x|\big)\in A(\mathbb{R}^{d})$. Тогда при $d_{1}=\Big[\frac{d-1}{2}\Big]$
(целая часть)
\begin{equation*}
    f_{d}\in C^{d_{1}}(0,+\infty),\quad
    \lim\limits_{t\rightarrow+\infty}t^{p}f^{(p)}_{d}(t)=0\quad (0\leq p\leq d_{1}),
\end{equation*}
а если $d_{1}\geq1$ ($d\geq3$), то и
\begin{equation*}
    \lim\limits_{t\rightarrow+0}t^{p}f^{(p)}_{d}(t)=0\quad (1\leq p\leq d_{1}).
\end{equation*}}

\emph{Кроме того, при $t>0$ сходится интеграл
\begin{equation*}
    \int\limits_{\rightarrow0}^{t}\Big(f^{(d_{1})}_{d}(t+u)-f^{(d_{1})}_{d}(t-u)\Big)\frac{du}{u}.
\end{equation*}}

Это необходимые условия [\ref{Trigub_1980}]. Там же есть и достаточные условия. См.
также замечание 3 о выпуклых функциях.

Давно известно, что для любой функции $f_{1}$ при любом $d\geq2$
\begin{equation*}
    \int\limits_{0}^{1}\big(1-u^{2}\big)^{\frac{d-3}{2}}f_{1}(ut)du=f_{d}(t).
\end{equation*}

Кроме того, при любой функции $f_{d}$ и нечетном $d\geq3$
$\Big(d_{1}=\frac{d-1}{2}\Big)$
\begin{equation*}
    \Big\{t^{d_{1}-\frac{1}{2}}f_{d}\big(\sqrt t\big)\Big\}^{(d_{1})}=f_{1}\big(\sqrt
    t\big)\quad (t\geq0)
\end{equation*}
(см. [\ref{Trigub_Bel}]), а при четном $d$ появляются производные полуцелого порядка
[\ref{Trigub_2010}]:
\begin{equation*}
    \sqrt t\int\limits_{0}^{t}\frac{d^{n+1}}{du^{n+1}}\Big\{u^{d_{1}}f_{d}\big(\sqrt
    u\big)\Big\}\frac{du}{\sqrt{t-u}}=f_{1}\big(\sqrt t\big)\quad \Big(t\geq0,\
    n=d_{1}=\frac{d}{2}-1\Big).
\end{equation*}

При этом взаимно однозначном соответствии из того, что $f_{d}=0$ на интервале, следует,
что и $f_{1}=0$ там же. А если $f_{1}(t)=0$ при $t\geq1$, то для того чтобы и $f_{d}=0$
при $t\geq1$ и нечетном $d\geq3$ необходимо и достаточно, чтобы
\begin{equation*}
\int\limits_{0}^{1}t^{2k}f_{1}(t)dt=0\quad \Big(0\leq k\leq\frac{d-3}{2}\Big),
\end{equation*}
так как в этом случае при $t\geq1$
\begin{equation*}
\int\limits_{0}^{1}f_{1}(ut)\big(1-u^{2}\big)^{\frac{d-3}{2}}du=\frac{1}{t}\int\limits_{0}^{1}f_{1}(u)\Big(1-\frac{u^{2}}{t^{2}}\Big)^{\frac{d-3}{2}}du.
\end{equation*}

Как хорошо известно, функция $f_{d}$ характеризуется следующим представлением
\begin{equation*}
    f_{d}(t)=\int\limits_{0}^{\infty}j_{\frac{d}{2}-1}(ut)g(u)du,\quad g\in
    L_{1}(\mathbb{R}_{+}),
\end{equation*}
где бесселева функция
\begin{equation*}
    j_{\lambda}(t)=\int\limits_{0}^{1}\big(1-u^{2}\big)^{\lambda-\frac{1}{2}}\cos
    utdu=c(\lambda)\frac{J_{\lambda}(t)}{t^{\lambda}}\quad
    \Big(\lambda>-\frac{1}{2}\Big)
\end{equation*}
и $j_{-\frac{1}{2}}(t)=\cos t$.

А если $f_{d}\big(|x|\big)\in A^{+}(\mathbb{R}^{d})$, то дополнительно $g(u)\geq0$ на
$\mathbb{R}_{+}$. Поэтому, если $f_{1}\big(|x|\big)\in A^{+}(\mathbb{R}^{1})$, то и
соответствующая функция $f_{d}\big(|x|\big)\in A^{+}(\mathbb{R}^{d})$, т.к.
\begin{equation*}
    \begin{split}
       & f_{d}(t)=\int\limits_{0}^{\infty}f_{1}(ut)\big(1-u^{2}\big)^{\frac{d-3}{2}}du=\int\limits_{0}^{1}\big(1-u^{2}\big)^{\frac{d-3}{2}}du\int\limits_{0}^{\infty}g(v)\cos uvtdv= \\
       &=\int\limits_{0}^{\infty}g(v)dv\int\limits_{0}^{1}\big(1-u^{2}\big)^{\frac{d-3}{2}}\cos
       utdu=\int\limits_{0}^{\infty}j_{\frac{d}{2}-1}(vt)g(v)dv.
    \end{split}
\end{equation*}

Противоположное утверждение (от $f_{d}$ к $f_{1}$) не верно, т.к. если $f_{1}\not\equiv
0$ и положительно определенная, то $f_{1}(0)>0$.

Для доказательства положительности $\widehat{f}$ используют два метода.

I. Неравенства типа
\begin{equation*}
\int\limits_{0}^{t}(t-u)^{\mu}J_{\nu}(u)u^{\nu}du\geq0\quad (t>0)
\end{equation*}
(ссылки на авторов таких неравенств см., напр., в [\ref{Zastavn_Trigub}]).

II. Если преобразование Лапласа при $t>0$ и некотором $\mu$
\begin{equation*}
    F_{1}(t)=\int\limits_{0}^{\infty}e^{-tu}u^{\mu}\widehat{f}_{1}(u)du\geq0
\end{equation*}
и интеграл является вполне монотонной функцией, т.е. $F_{1}\in C^{\infty}(0,+\infty)$ и
$(-1)^{k}F_{1}^{(k)}(t)\geq0$ для $t>0$ и любых целых $k\geq0$, то $\widehat{f}_{1}(t)\geq0$
$(t>0)$. См., напр., [\ref{Zastavn_Trigub}].

В этом случае получается и двусторонняя оценка для $\widehat{f}_{1}(t)$, что является важным
для оценки порядка приближения линейными комбинациями сдвигов одной функции (см.
[\ref{Zastavn_Trigub}], [\ref{Buhmanm}]).

\textbf{\emph{Замечание 5.}} (о других нормах в $\mathbb{R}^{d}$)

\emph{Как известно, преобразование Фурье --- оператор, коммутирующий с вращением
пространства $\mathbb{R}^{d}$ вокруг нуля. А после поворота квадрата на угол
$\frac{\pi}{4}$ получаем квадрат $|x_{1}|+|x_{2}|\leq\sqrt 2$.}

Так что теоремы \ref{Th_3} -- \ref{Th_7} можно переформулировать для функций вида
$f_{0}\big(|x_{1}|+|x_{2}|\big)$.

К настоящему времени известно уже много результатов о положительной определенности
функций, зависящих от одной какой--либо нормы в $\mathbb{R}^{d}$ при разных размерностях
$d$. См. [\ref{Cambanis_Keener_Simons}] и [\ref{Zastavn}] с большим списком литературы.

\textbf{\emph{E--mail address}}: roald.trigub@gmail.com

\end{document}